\documentclass{siamltex}

\usepackage[
  hyperindex  = {true},
  colorlinks  = {true},
  linkcolor   = {blue},
  citecolor   = {blue}
]{hyperref}
\usepackage{multirow}
\usepackage{color, colortbl}
\usepackage{caption}
\usepackage{subcaption}
\usepackage{latexsym,graphicx, amssymb, amsfonts, setspace, geometry, amsmath}
\usepackage{pgfplots}
\pgfplotsset{compat=1.8}
\usepgfplotslibrary{statistics}
\usepackage{grffile}
\usepackage{pgfplotstable}
\usepackage{slashbox}
\usepackage{mathtools}
\usepackage[utf8]{inputenc}
\usepackage[english]{babel}
\usepackage{float}
\usepackage{esvect}
\usepackage{algorithm,algorithmic}
\usepackage{enumitem}
\usepackage{hhline}
\usepackage{varwidth}
\usepackage{tikz}
\usetikzlibrary{backgrounds,patterns,calc}
\usepackage{xparse}
\usepackage{tabularx}
\usepackage{placeins}

\pgfplotsset{compat=newest}

\makeatletter
\newcommand{\doublewidetilde}[1]{{%
  \mathpalette\double@widetilde{#1}%
}}
\newcommand{\double@widetilde}[2]{%
  \sbox\z@{$\m@th#1\widetilde{#2}$}%
  \ht\z@=.9\ht\z@
  \widetilde{\box\z@}%
}
\makeatother

\newcommand{\xb}{\mathbf{x}}
\newcommand{\ub}{\mathbf{u}}

\newcommand{\db}{\mathbf{d}}

\newcommand{\Ib}{\mathbf{I}}

\newcommand{\tran}{^{\top\kern-\scriptspace}}

\newcommand{\be}{\begin{equation}}
\newcommand{\ee}{\end{equation}}
\newcommand{\ba}{\begin{aligned}}
\newcommand{\ea}{\end{aligned}}
\newcommand{\bea}{\begin{eqnarray}}
\newcommand{\eea}{\end{eqnarray}}

\newcommand\nc{\newcommand}

\newcommand{\bE}{\mathbf{E}}
\newcommand{\bH}{\mathbf{H}}

\nc\ex{E_x}
\nc\ey{E_y}
\nc\ez{E_z}

\nc\hx{H_x}
\nc\hy{H_y}
\nc\hz{H_z}

\nc\px[1]{\frac{\partial #1}{\partial x}}
\nc\py[1]{\frac{\partial #1}{\partial y}}

\nc\br{{\mathbf{r}}}
\nc\curl{\nabla\times}
\nc\dive{\nabla\cdot}
\nc\pt{\frac{\partial}{\partial t}}
\nc\pet{\frac{\partial \bE}{\partial t}}
\nc\pht{\frac{\partial \bH}{\partial t}}

\newcommand{\bx}{\mathbf{x}}
\newcommand{\by}{\mathbf{y}}

\DeclareMathOperator*{\argmin}{arg\,min}

\definecolor{Gray}{gray}{0.9}
\newcolumntype{g}{>{\columncolor{Gray}}c}

\pgfmathdeclarefunction{gauss}{3}{%
  \pgfmathparse{1/(#3*sqrt(2*pi))*exp(-((#1-#2)^2)/(2*#3^2))}%
}

\title{Inverse scattering reconstruction of a three dimensional sound-soft axis-symmetric impenetrable object}

\author{Carlos Borges\thanks{Department of Mathematics, University of Central Florida, Orlando, FL, USA. \textit{Email: Carlos.Borges@ucf.edu}}
\and
Jun Lai\thanks{School of Mathematical Sciences, Zhejiang University, Hangzhou, Zhejiang, China. \textit{Email: laijun6@zju.edu.cn}}
}

\begin{document}

\maketitle

\begin{abstract}
In this work, we consider the problem of reconstructing the shape of a three dimensional impenetrable sound-soft axis-symmetric obstacle from measurements of the scattered field at multiple frequencies. This problem has important applications in locating and identifying obstacles with axial symmetry in general, such as, land mines. We present a two-part framework for recovering the shape of the obstacle. In part 1, we introduce an algorithm to find the axis of symmetry of the obstacle by making use of the far field pattern. In part 2, we recover the shape of the obstacle by applying the recursive linearization algorithm (RLA) with multifrequency measurements of the scattered field. In the RLA, a sequence of inverse scattering problems using increasing single frequency measurements are solved. Each of those problems is ill-posed and nonlinear. The ill-posedness is treated by using a band-limited representation for the shape of the obstacle, while the nonlinearity is dealt with by applying the damped Gauss-Newton method. When using the RLA, a large number of forward scattering problems must be solved. Hence, it is paramount to have an efficient and accurate forward problem solver. For the forward problem, we  apply separation of variables in the azimuthal coordinate and Fourier decompose the resulting problem, leaving us with a sequence of decoupled simpler forward scattering problems to solve. Numerical examples for the inverse problem are presented to show the feasibility of our two-part framework in different scenarios, particularly for objects with non-smooth boundaries.
\end{abstract}
\section{Introduction}\label{s:intro}
\begin{figure}[h]
\begin{subfigure}[t]{0.48\textwidth}
\center
\includegraphics[width=1\textwidth]{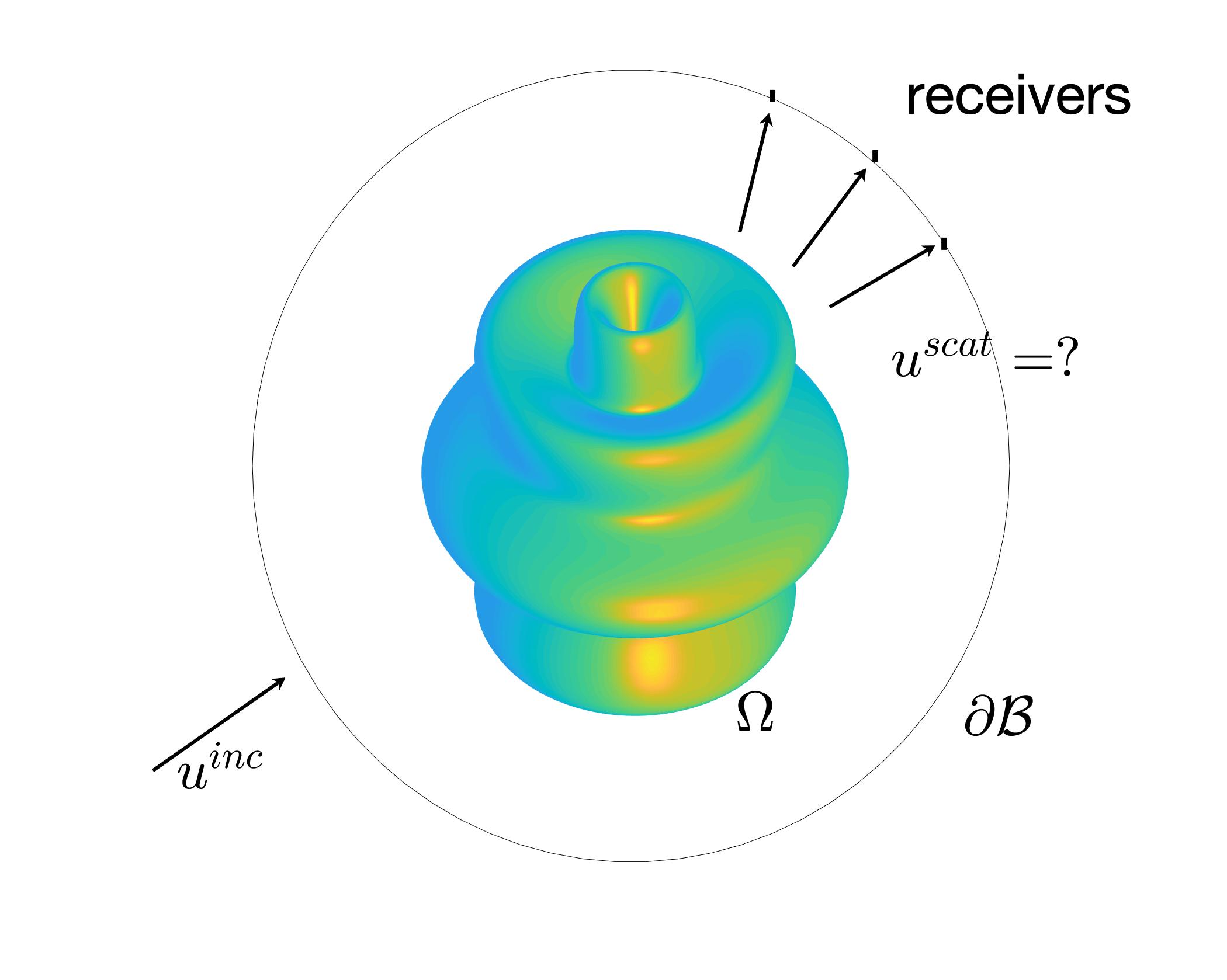}
\caption{Forward scattering problem}\label{fig:fwd_prob}
\end{subfigure}
\begin{subfigure}[t]{0.48\textwidth}
\center
\includegraphics[width=1\textwidth]{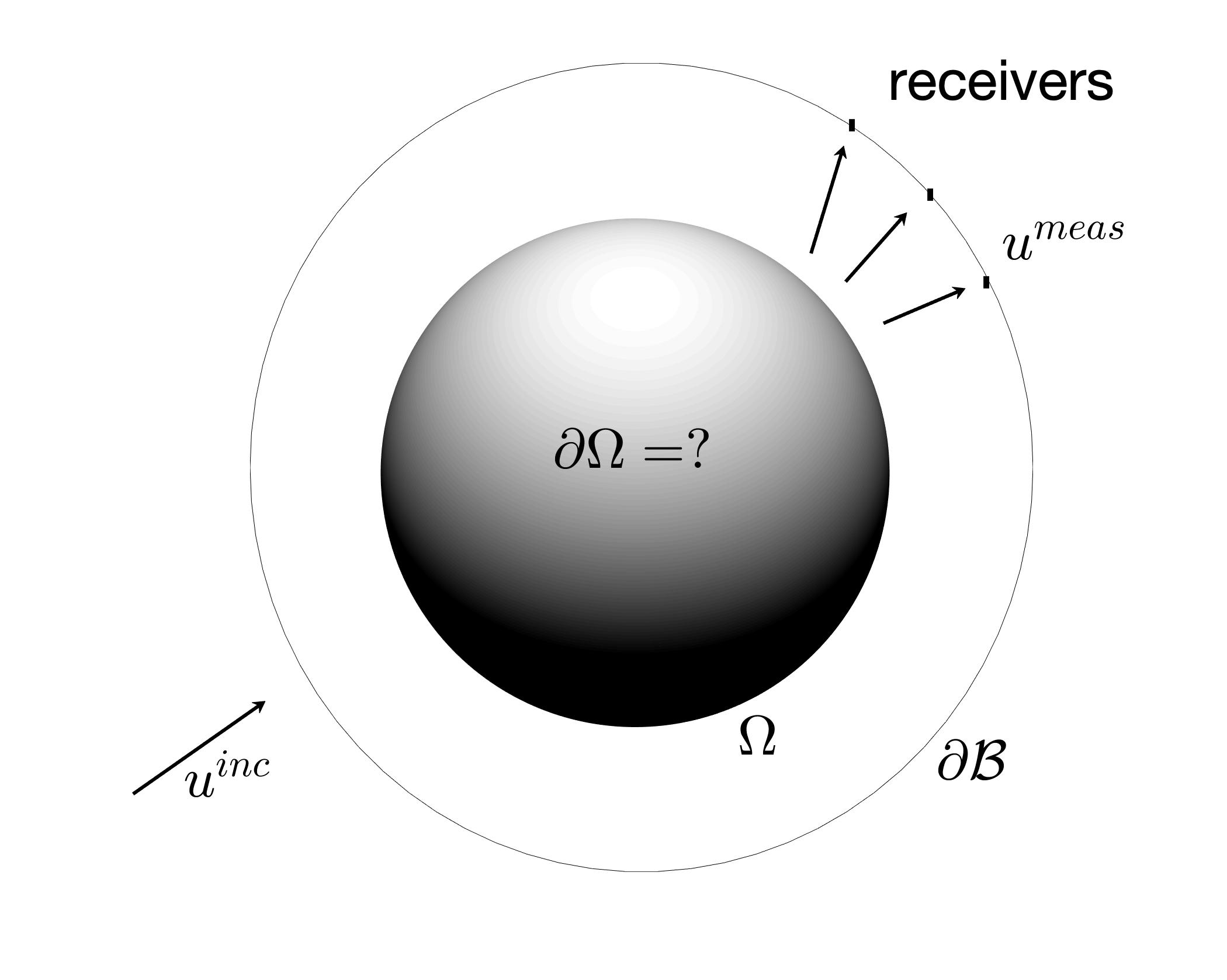}
\caption{Inverse scattering problem}\label{fig:inv_prob}
\end{subfigure}
\caption{Scattering from an axis-symmetric sound-soft obstacle. In the {\em forward scattering problem}, the boundary $\partial \Omega$ representing the shape of the obstacle $\Omega$ is known and we want to evaluate the scattered field $u^\emph{scat}$ given the incident field $u^\emph{inc}$, as shown in Figure \ref{fig:fwd_prob}. In the {\em inverse scattering problem}, the shape of the obstacle $\partial \Omega$ is unknown and we want to determine it using measurements of the scattered field $u^\emph{meas}$ at the receivers located on $\partial\mathcal{B}$, as shown in Figure \ref{fig:inv_prob}.}\label{fig:problems_intro}
\end{figure}

There are a large amount of important applications of inverse scattering, such as medical imaging \cite{Hoskins,kuchment2014radon,
nashed2002inverse,scherzer2010handbook,simonetti2008inverse}, nondestructive testing \cite{collins1995nondestructive,engl2012inverse,langenberg1993imaging}, remote sensing \cite{Ustinov2014}, ocean acoustics \cite{0266-5611-10-5-003}, geophysics \cite{zhdanov2002geophysical,snieder1999inverse,Eldad_04}, sonar and radar \cite{cheney2009fundamentals,Colton}, and many others. Among those applications, the recovery of the shape of axis-symmetric or nearly axis-symmetric obstacles and cavities plays a very important role in practice, as for instance, the identification and classification of locations and types of different missiles and  mines. In this paper, we consider the forward and inverse scattering problems in three dimensions for an axis-symmetric sound-soft obstacle $\Omega$, as described in Figure \ref{fig:problems_intro}. 

%

We define the forward scattering operator for this problem as the operator $\mathcal{F}_{k,\db}:\partial \Omega \rightarrow \mathbb{C}^{M}$, such that
\begin{equation} \label{eq:fwd_op}
\mathcal{F}_{k,\db}(\partial \Omega)=\ub_{k,\db}^\emph{meas}\rvert_{\partial \mathcal{B}}
\end{equation}
where $\partial \Omega$ is the boundary of the obstacle, and  $\ub_{k,\db}^\emph{meas}\vert_{\partial \mathcal{B}}$ is a vector in $\mathbb{C}^{M}$ with coordinates being the measurements of the scattered field $u_{k,\db}^\emph{scat}$ at $M$ receptors located on a surface $\partial \mathcal{B}$. The scattered field $u_{k,\db}^\emph{scat}$ is generated by the incidence of a plane wave $u_{k,\db}^\emph{inc}$ with incident direction $\db$ and wavenumber $k$ and can be obtained by solving the Helmholtz equation
%
%
\begin{alignat}{4}\label{eq:helm}
\Delta u_{k,\db} + k^2 u_{k,\db} &= 0, &\quad &\text{in}\quad \mathbb{R}^3\setminus\Omega,\\\nonumber
                       u_{k,\db} &= 0, &\quad &\text{on}\quad \partial \Omega, \nonumber
\end{alignat}
where $u_{k,\db} = u_{k,\db}^\emph{scat}+u_{k,\db}^\emph{inc}$ is the total field. The scattered wave $u_{k,\db}^\emph{scat}$ also satisfies the Sommerfeld radiation condition
\begin{equation*}
\lim_{r\rightarrow \infty}r\bigg(\frac{\partial u_{k,\db}^\emph{scat}}{\partial \nu}-iku_{k,\db}^\emph{scat}\bigg) = 0, \quad r=|x|,
\end{equation*}
where $\nu$ is the exterior unit normal of $\Omega$. We assume the wavenumber $k$ satisfies $\Re(k)>0$ and $\Im(k)\ge 0$.   To solve the Helmholtz equation, one can apply potential theory to obtain the integral equation formulation of \eqref{eq:helm} as described in \cite{Colton}. When the obstacle $\Omega$ is arbitrary, the evaluation of the integral operators defined on a surface in three dimensions requires a very costly treatment of the quadratures and discretization of the boundary. Moreover, solving the integral equation requires specific schemes, such as fast multipole method together with a Krylov subspace iterative method, like GMRES \cite{Greengard87,carrier1988fast}.

On the other hand, when the obstacle is axis-symmetric, the forward solver can be greatly simplified. Since $\partial \Omega$ is obtained by rotating a curve $\gamma$ around the axis of symmetry, in a slight abuse of notation, we rewrite the forward operator \eqref{eq:fwd_op} as
\begin{equation} \label{eq:fwd_op_curve}
\mathcal{F}_{k,\db}(\gamma)=\ub_{k,\db}^\emph{meas}\vert_{\partial \mathcal{B}},
\end{equation}
where $\gamma:[0,1]\rightarrow\mathbb{R}^2$ is an open simple curve with both end points on the axis of symmetry of the obstacle. Due to the symmetry, solving the forward problem can be accelerated by applying separation of variables in the azimuthal angle and Fourier decomposing the resulting integral equation \cite{Gedney1990, HAO2015304, HELSING2014686, LIU2016226, young2012high}. The original integral equation turns into a sequence of uncoupled integral equations, one for each Fourier mode, and the integral operators for these equations are defined along the curve $\gamma$ only. In this case, the quadrature scheme for each integral equation on $\gamma$ is much easier to implement and the system of equations is much  cheaper to solve.

In this paper, we are interested in reconstructing the shape of an axis-symmetric obstacle given measurements of the scattered field at the receivers from one or more incident waves. We propose a two-part framework for recovering the shape of the obstacle. In part 1 of the framework, we introduce a method that uses the full-aperture data from two incoming incident waves to obtain the axis of symmetry of the obstacle by looking into the symmetry of the far field along a circle in a plane. In part 2 of the framework, we apply the recursive linearization algorithm (RLA) \cite{bao2005inverse, bao2007inverse, bao2010error, bao2015inverse, Chen, chen1995recursive, sini2012inverse, sini2013convergence, Borges2015, Borges2017} to recover the curve $\gamma$. In doing this, a sequence of single frequency inverse problems of the form
\begin{equation}
\tilde{\gamma}=\arg \min_\gamma \| \ub_k^{\emph{meas}} -\mathcal{F}_k(\gamma)\|, \label{eq:single_freq_inv}
\end{equation}
is solved, where 
\begin{equation*}
\ub_k^{\emph{meas}}=\left[u_{k,\db_1}^{\emph{meas}}; \cdots;u_{k,\db_{N_d}}^{\emph{meas}}\right] \quad \text{and} \quad \mathcal{F}_k=\left[\mathcal{F}_{k,\db_1}; \cdots;\mathcal{F}_{k,\db_{N_d}}\right].
\end{equation*}
The RLA works as a continuation method in the wavenumber parameter, where we use the reconstruction from the previous frequency as the initial guess for the next one. As the single frequency inverse problem \eqref{eq:single_freq_inv} is highly nonlinear and ill-posed \cite{Colton},  we propose a damped Gauss-Newton method combined with a band-limited regularization of the curve $\gamma$ as in \cite{Borges2015} to overcome the difficulties. In the end, since we fully make use of the symmetry, our algorithm is extremely efficient and can accurately locate and reconstruct the unknown object, even with nonsmooth boundary.  
  
{\bf Related work:} We refer readers to \cite{Gedney1990, HAO2015304, HELSING2014686, LIU2016226, young2012high} for the forward acoustic and electromagnetic scattering problems for axis-symmetric obstacles. The inverse scattering problem for three dimensional obstacles was studied in \cite{Farhat2002,harbrecht2007fast,gutman1993regularized,gutman1993two,gutman1994iterative}. The time domain inverse scattering problem for three dimensional obstacles was studied in \cite{beilina2015globally,thanh2015imaging}. Readers are referred to \cite{bao2005inverse, bao2007inverse, bao2010error, Borges2015, Borges2017, CHAILLAT20124403, chen1995recursive, Chen, sini2012inverse, sini2013convergence} for the inverse scattering problem for two and three dimensions using multiple frequency data. In particular, a complete review on inverse scattering problems based on multiple frequency data was given in \cite{bao2015inverse}. Recently, authors in \cite{shin2020} proposed an algorithm to determine the two dimensional radially symmetric potential from single frequency near-field scattering data.  However, we are not aware of any previous work on the three dimensional inverse obstacle problem using multiple frequency data when the obstacle has an axis of symmetry. 

{\bf Contributions:} The contributions of this paper are summarized as follows:
\begin{itemize}
	\item We obtain a uniqueness result for the inverse scattering of an axis-symmetric object with single frequency data by plane wave incidence. 
	\item We propose a novel algorithm to determine the orientation and location of the axis of symmetry of the unknown object based on single frequency data.
	\item We apply the recursive linearization algorithm with multifrequency data and band-limited representation to reconstruct the generating curve of the axis-symmetric obstacle.
\end{itemize}

{\bf Notation:} We present the most common symbols used in this paper in Table \ref{table:symbol}. 
\begin{table}
\caption{List of main symbols used in this article.}\label{table:symbol}
{\small
\begin{center}
\begin{tabular}{ll}
\hline
Symbol & Description \\
\hline\hline
$\Omega$                & Closed set representing the sound-soft impenetrable obstacle\\
$\partial \Omega$     & Boundary of the obstacle $\Omega$\\
$\gamma$                & Parametrization of the curve used to generate the axis-symmetric obstacle \\
$\db$                        & Incident direction of plane wave $u^{\emph{inc}}_{k,\db}$ ($\|\db\|=1$) \\
$k$                           & Wavenumber (or frequency) of the incident plane wave \\
$u_{k,\db}^{\emph{inc}}$  & Incident plane wave with wavenumber $k$ and incident direction $\db$\\
$u_{k,\db}^{\emph{scat}}$ & Scattered field off of the obstacle $\partial \Omega$ generated by $u_{k,\db}^{\emph{inc}}$ \\
${\bf u}_{k,\db}^{\emph{meas}}$  & Vector with coordinates being $u_{k,\db}^{\emph{scat}}$ measured at the receivers\\
$u_{k,\db}^{\infty}$  & Far-field pattern of the scattered field $u_{k,\db}^{\emph{scat}}$\\
$N$                           & Number of discretization points at the boundary  $\gamma$ of the obstacle\\
$M$				& Number of receptors \\
$N_d$                       & Number of incident waves \\
$\mathcal{F}_{k,\db}$  & Forward scattering operator mapping $\partial\Omega$ to $u_{k,\db}^\emph{scat}$ (for given $u_{k,\db}^{\emph{inc}}$) \\
$\mathcal{F}_k$  & Forward scattering operator at wavenumber $k$ for $N_d$ directions ($\left[\mathcal{F}_{k,\db_1}; \cdots;\mathcal{F}_{k,\db_{N_d}}\right]$)\\
$\partial_\gamma \mathcal{F}_k$ & Frech\'{e}t derivative of $\mathcal{F}_k$ with respect to $\gamma$ \\
$\mathcal{S}$           & Single layer potential  \\
$\mathcal{D}$           & Double layer potential \\
$\mathcal{I}$             & Identity operator \\
$G^{k}$                     &  Free space Green's function for the three dimensional Helmholtz equation \\
$G_m^{k}$                & Modal Green's function for the $m^{th}$ mode\\
$\mathcal{S}_m$  & Modal single layer potential for the $m^{th}$ mode \\
$\mathcal{D}_m$  & Modal double layer potential for the $m^{th}$ mode\\
$\bf S_m$           & $N\times N$ matrix for the discretization of the $m^{th}$ modal single layer operator \\
$\bf D_m$           & $N\times N$ matrix for the discretization of the $m^{th}$ modal double layer operator \\
\hline 
\end{tabular}
\end{center}
}
\end{table}

{\bf Article Outline:} In Section \ref{s:fwd_problem}, we introduce the fast solver for the  forward scattering problem of a three dimensional axis-symmetric obstacle. In Section \ref{s:inv_problem}, we show the uniqueness result for the inverse axis-symmetric obstacle problem and propose a two-step framework to reconstruct the shape of the unknown obstacle. In Section \ref{s:num_res}, numerical examples are presented to illustrate different characteristics of the method. Concluding remarks are made in Section \ref{s:conclusions}.

\section{Forward Scattering Problem}\label{s:fwd_problem}
To evaluate the forward scattering operator, we must solve the problem \eqref{eq:helm} for $u_{k,\db}^\emph{scat}$, given an incident wave $u_{k,\db}^\emph{inc}$. As we are considering the forward problem with a fixed wavenumber $k$ and direction $\db$ in this section and the next one, to ease the notation and when there is no confusion, we will drop the indices for $k$ and $\db$ for the fields unless it is otherwise stated. We represent the scattered field using layer potentials. First, we define, respectively,  the single and double layer potentials for $\bx\in \mathbb{R}^3\backslash\Omega$ as
\begin{eqnarray*}
\mathcal{S}\mu(\bx) = \int_{\partial \Omega}G^k(\bx,\by)\mu(\by)ds(\by),\quad \text{and}\quad
\mathcal{D}\mu(\bx) = \int_{\partial \Omega}\frac{\partial G^k(\bx,\by)}{\partial \nu(\by)}\mu(\by)ds(\by),
\end{eqnarray*}
where $\nu$ is the exterior unit normal to the boundary of the obstacle and $G^k(\bx,\by)$ is the free space Green's function of Helmholtz equation, i.e.
\begin{eqnarray*}
G^k(\bx,\by) = \frac{e^{ik|\bx-\by|}}{|\bx-\by|}.
\end{eqnarray*}
To avoid resonances, we chose to represent the scattered field using a combined layer potential approach and write
\begin{eqnarray}\label{eq:comb_lay}
u^\emph{scat}(\bx) = \left(\mathcal{D}+ik\mathcal{S}\right)\mu(\bx).
\end{eqnarray}

Using \eqref{eq:comb_lay} and the sound-soft boundary condition with the jump properties from Theorem 3.1 in \cite{Colton}, we obtain a uniquely solvable equation for any $k$ with $\Re(k)>0$ and $\Im(k)>0$,  
\begin{eqnarray} \label{eq:comb_layer}
\left(\frac{1}{2}\mathcal{I}+\mathcal{D}+ik\mathcal{S}\right)\mu(\bx) = -u^\emph{inc}(\bx) 
\end{eqnarray}
for $\bx \in \partial \Omega$.

Next, we make use of the axis-symmetry of $\partial \Omega$ and rewrite the density function $\mu(\bx)$ as
\begin{eqnarray} \label{eq:potmodal}
\mu(r,\theta, z) = \sum_{m = -\infty}^{ \infty}\mu_m(r, z) e^{im\theta}, 
\end{eqnarray} 
and the incident plane wave as
\begin{equation} \label{eq:incmodal}
u^\emph{inc}(r,\theta,z) = \sum_{m = -\infty}^{ \infty}u^\emph{inc}_m(r, z) e^{im\theta},
\end{equation}
where $(r,\theta, z)$ are the cylindrical coordinates of $\bx\in\partial \Omega$. Using \eqref{eq:potmodal} and \eqref{eq:incmodal} in Equation \eqref{eq:comb_layer}, we obtain a sequence of line integral equations
\begin{equation}\label{eq:modal_eqs}
\left(\frac{1}{2}\mathcal{I}+\mathcal{D}_m+ik\mathcal{S}_m\right)\mu_m =  -u^\emph{inc}_m, \mbox{ for } m \in \mathbb{Z}
\end{equation}
with 
\begin{eqnarray}
\mathcal{S}_m \mu_m &=& \int_{\gamma}G_m^k(r,z,r',z')\mu_m(r', z')r'ds(r',z'), \label{eq:Skm}  \\
\mathcal{D}_m \mu_m &=& \int_{\gamma}\frac{\partial G_m^k(r,z,r',z')}{\partial \nu(r',z')} \mu_m(r', z')r'ds(r',z'), \label{eq:Dkm}
\end{eqnarray}
where  $G_m^k(r,z,r',z')$ are the modal Green's functions given by 
\begin{equation}\label{eq:modeS}
G_m^k(r,z,r',z') = \int_{0}^{2\pi}G^k(r,0,z,r',\theta',z')d\theta'= \int_0^{2\pi}\frac{e^{ik\rho}}{\rho} e^{-im\theta'}d\theta',
\end{equation}
with $\rho = \sqrt{r^2+r'^2+(z-z')^2-2rr'\cos\theta'}$ and 
\begin{equation} \label{eq:modeD}
\frac{\partial G_m^k(r,z,r',z')}{\partial \nu} = \int_0^{2\pi}(ik\rho-1)(\nu_{r'}(r\cos\theta-r')+\nu_{\theta'}(z-z')) \frac{e^{ik\rho}}{\rho^3} e^{-im\theta'}d\theta'.
\end{equation}
It is worth mentioning that all the integral equations in \eqref{eq:modal_eqs} are decoupled from each other, which greatly simplifies the computation of the forward problem.

To solve each of the integral equations in \eqref{eq:modal_eqs}, we must evaluate the line integrals in \eqref{eq:modeS} and \eqref{eq:modeD}. Unfortunately, the modal Green's functions are not in closed form and the kernels in these integral operators have strong singularities. To accelerate the computation, we apply an FFT based algorithm \cite{LAI20171, LAI2019} with recursive formulas to efficiently evaluate the modal Green's functions. In order to discretize the singular integral \eqref{eq:Skm} and \eqref{eq:Dkm} to high order,  we divide the curve $\gamma$ into a set of panels, such that each panel has at least 12 points per wavelength. Next, we discretize each panel using 16 Gauss-Legendre nodes, and use the $16^{th}$-order generalized Gaussian quadrature from \cite{bremer2010nonlinear} to apply the Nystr\"om method in each of the equations \eqref{eq:modal_eqs} to handle the singularity in the modal Green's functions. In the end, for each mode, we obtain an $N\times N$ system of linear equations given by
\begin{equation*} 
(\Ib+{\bf D_m} +ik {\bf S_m}) {\boldsymbol \mu}_m = {\bf -u^\emph{inc}_m}
\end{equation*}
where $\Ib$ is the $N\times N$ identity matrix, ${\bf D_m}$ and ${\bf S_m}$ are the $ N\times N$ matrices obtained by the discretization of the potentials $\mathcal{D}_m$ and $\mathcal{S}_m$, ${\boldsymbol \mu}_m$ and ${\bf -u^\emph{inc}_m}$ are the $\mathbb{C}^N$ vectors with coordinates being the values of the density $\mu_m$ and the function $-u^\emph{inc}_m$, respectively, at the discretization points on $\gamma$.

For the computational complexity, if the size of the obstacle is $\mathcal{O}(1)$, we have that $N=\mathcal{O}(k)$ and the system is small enough to be efficiently solved using Gaussian Elimination on $\mathcal{O}(k^3)$ operations. The number of modes that need to be calculated to resolve the scattered field, which is also the number of linear systems that need to be solved,  is $\mathcal{O}(k)$. The total work to calculate the scattered field for a single incident wave is $\mathcal{O}(k^4)$. If the scattered field needs to be calculated for $N_d$ incident waves, the total work becomes $\mathcal{O}(k^4+N_dk^3)$, where the first term comes from the calculation of the inverse matrices for all modes and the second refers to the application of those inverse matrices in the $N_d$ incoming waves. It is much more efficient than a general 3D forward solver which usually has complexity on the order of $\mathcal{O}(k^6)$.

\section{Inverse Scattering Problem}\label{s:inv_problem}
A large family of scattering objects, in practice, can be represented by shapes obtained by rotating a curve along an axis. This representation, even though it has its limitations, covers several important applications, such as nano particles, industrial machinery parts and missiles. Suppose that $\gamma:[0,1]\rightarrow\mathbb{R}^2$ is a parametrization of the curve rotated along the axis to generate the boundary of the obstacle $\Omega$. Given the forward problem \eqref{eq:fwd_op_curve}, we are interested in the following inverse problem:

{\it {\bf Inverse Obstacle Problem (Axis-symmetric case):}  Given the measurements ${\bf u}_{k_i,\db_j}^\emph{meas}$ of the scattered field $u_{k_i,\db_j}^\emph{scat}$ of an unknown impenetrable axis-symmetric obstacle $\Omega$ for some known collection of incident waves $u_{k_i,\db_j}^\emph{inc}=e^{ik_i\xb\cdot\db_j}$, $i=1,\ldots,N_k$, $j=1,\ldots,N_d$, obtain a reconstruction of the shape of $\Omega$.}

For a general three dimensional obstacle, given the scattered field $u^\emph{scat}$ of the obstacle generated by the scattering of an incident plane wave $u^\emph{inc}=e^{ik\bx\cdot\db}$, one cannot expect to uniquely recover the shape of the obstacle\cite{Colton}. However, for the case of an axis-symmetric object, if we assume the axis of symmetry is fixed, then measurements based on one incident plane wave are enough to determine the shape of the obstacle, as shown in Theorem \ref{thm:uniq_axis}. 

\begin{theorem} \label{thm:uniq_axis}
	 Assume  $\Omega_1$ and $\Omega_2$ are two axis-symmetric scatterers with the same axis such that the scattering fields $u_{1}^{\emph{scat}}$ and $u_{2}^{\emph{scat}}$ on $\partial \mathcal{B}$ coincide for one incident plane wave $u^{\emph{inc}}=e^{ik\bx\cdot\db}$ with $\db\ne(0,0,\pm1)$. Then $ \Omega_1 = \Omega_2$.
\end{theorem}
\begin{proof}
	Assume $\Omega_1\ne \Omega_2$. Since the scattered field on $\partial \mathcal{B}$ uniquely determines the far field and the far field uniquely determines the the scattered field outside the region $\Omega_1\cup\Omega_2$, we have $u_{1}^\emph{scat}(\bx)=u_{2}^\emph{scat}(\bx)$ for $\bx\in \mathbb{R}^3\backslash(\Omega_1\cup\Omega_2)$. Without loss of generality, we may assume $D=(\mathbb{R}^3\backslash\Omega_2)\cap\Omega_1$ is nonempty. Inside $D$, we have $u = u_{2}^\emph{scat}+u^\emph{inc}$ well defined and it satisfies the zero boundary condition on $\partial D$. Therefore, $u$ is a Dirichlet eigenfunction for the negative Laplacian in the domain $D$ with eigenvalue $k^2$. Next, we will show that this implies there exists infinitely many eigenfunctions for the same eigenvalue $k^2$. 
	
	Since both scatterers are axis-symmetric and share the same axis, $D$ is also axis-symmetric. We apply the Fourier decomposition along the azimuthal direction to the incident field $u^\emph{inc}$. Let $\db=(d_1,d_2,d_3)$ and $\bx=(r\cos\theta,r\sin\theta,z)$. According to the Jacobi-Anger formula \cite{Colton}, the plane wave has the expansion
	\begin{eqnarray}
	u^\emph{inc}(r,\theta,z)=e^{ik(d_1r\cos\theta+d_2r\sin\theta+d_3z)} = \sum_{m=-\infty}^{\infty}i^mJ_m(k\rho)e^{im(\theta-\phi)}e^{ikd_3z} \nonumber
	\end{eqnarray}
where $\phi = \arctan(d_1/d_2)$, $\rho = r\sqrt{d_1^2+d_2^2}$ and $J_m$ is the Bessel function of order $m$. In other words, the $m$-th mode of $u^\emph{inc}$ is
\begin{equation*}
u_{m}^\emph{inc}(r,\theta,z)=i^mJ_m(k\rho)e^{im(\theta-\phi)}e^{ikd_3z}.
\end{equation*}
For each $u_m^\emph{inc}$, the corresponding scattered field is given by the Fourier decomposition of $u^\emph{scat}$ along the azimuthal direction. We have that $u_{m} = u_{m}^\emph{inc}+u_{m}^\emph{scat} $ satisfies the zero boundary condition on $\partial D$. 

We show that $u_{m}$ is not identically zero in $D$ if $u_{m}^\emph{inc}$ is nonzero. If this is true, then $u_{m}^\emph{scat} = -u_{m}^\emph{inc}$ in $D$. By analyticity, $u_{m}^\emph{scat} = -u_{m}^\emph{inc}$ in $\mathbb{R}^3\backslash (\Omega_1\cup\Omega_2)$. This is a contradiction since $u_{m}^\emph{scat}$ satisfies the Sommerfeld radiation condition while $u_{m}^\emph{inc}$ does not. Thus, we obtain infinitely many linearly independent eigenfunctions in $D$ with eigenvalue $k^2$, which is a contradiction. Therefore, $\Omega_1=\Omega_2$. 
\end{proof}

From the conclusion of the previous theorem, we propose a two-part framework to find the shape of an axis-symmetric obstacle as follows:
\begin{itemize}
\item Part 1: find the axis of symmetry of the obstacle;
\item Part 2: recover the shape of the generating curve for the obstacle.
\end{itemize}

In part 1, we propose a procedure that will explore the symmetry of the obstacle to obtain the axis of symmetry. In particular, by inspecting the far field pattern of the scattered field of some incident waves with fixed frequency and different direction, we can determine the location of the axis of symmetry of the obstacle. In part 2, we apply the recursive linearization algorithm with band-limited representation to solve a sequence of inverse scattering problems using the wavenumber as a continuation parameter and obtain a high resolution reconstruction of the shape of the obstacle.

\subsection{Finding the axis of symmetry (Part 1)}
To be able to use the uniqueness result from Theorem \ref{thm:uniq_axis}, one must first find the axis of symmetry of the obstacle. Our algorithm to find the axis of symmetry is based on Theorems \ref{thm:rotation} and \ref{thm:location}.
\begin{theorem}\label{thm:rotation}
If the axis of symmetry of $\Omega$ is the $z$-axis, then for any $\phi \in [0,\pi]$, both the real and imaginary parts of the far field $u^{\infty}(\theta,\phi)$ of $\Omega$ by the incident wave $u^{inc}=e^{ik\db\cdot \bx }$ with $\db=(d_1,d_2,d_3)\ne(0,0,\pm 1)$ are symmetric with respect to $\theta = \theta_0$ and $\theta = 2\pi-\theta_0$, where $\theta_0$ satisfies $\cos(\theta_0)=d_1/\sqrt{d_1^2+d_2^2}$, and  $\sin(\theta_0)=d_2/\sqrt{d_1^2+d_2^2}$. Here $\theta\in[0,2\pi)$ is the  azimuthal angle in the $xy$-plane from the positive $x$-axis and $\phi\in[0,\pi]$ is the altitude angle from the positive $z$-axis.  If $\db=(0,0,\pm 1)$, then the far field is axis-symmetric with respect to the $z$-axis.
\end{theorem}
\begin{proof}
	For $\db\ne(0,0,\pm 1)$, without loss of generality, we may assume $d_2=0$, in which case we need to show that the far field is symmetric with respect to $\theta=0$ and $\theta = \pi$.  In fact, when $d_2=0$, we have $u^\emph{inc}=e^{ik(d_1r\cos\theta+d_3z)}=e^{ik(d_1r\cos(-\theta)+d_3z)}$, so  $u^\emph{inc}$ is symmetric with respect to $\theta=0$ and $\theta = \pi$ for any fixed $z$. On the other hand, the obstacle $\Omega$ is also symmetric with respect to the plane that is cut by $\theta=0$ and $\theta = \pi$. By the uniqueness theorem of the exterior problem, we have that the scattered field $u^\emph{scat}$ is also symmetric with respect to the plane where $\theta=0$ or $\theta = \pi$. The far field pattern can be obtained by using the equation
	\begin{eqnarray}\label{far_field}
	u^\infty(\theta,\phi)&=\frac{1}{4\pi}\int_{\partial \Omega}\bigg\{u^\emph{scat}(\by)\frac{\partial e^{-\mathrm{i}k \hat{\bx}\cdot \by}}{\partial \nu(\by)}-\frac{\partial u^\emph{scat}}{\partial \nu }(\by)e^{-\mathrm{i}k \hat{\bx}\cdot \by}\bigg\}\mathrm{d}s(\by),
	\end{eqnarray} 
	where $\hat{\bx}=(\cos\theta\sin\phi, \sin\theta\sin\phi, \cos\phi)$. Since all the components on the right hand side of equation \eqref{far_field} are symmetric with respect to $\theta=0$ or $\theta = \pi$, the far field must be symmetric with respect to $\theta=0$ or $\theta = \pi$, too.
	
	Similarly, when $\db=(0,0,\pm 1)$, both the incident wave and the obstacle are axis-symmetric with respect to the $z$-axis, so is the scattered field and the far field pattern.
	\end{proof}

Theorem \ref{thm:rotation} implies that the far field pattern is symmetric with respect to the axis of symmetry when the axis passes through the origin. For an axis-symmetric obstacle that is not centered at the origin, we have the following translation property for the far field. 
\begin{theorem}\label{thm:location}
	Let $u^{\infty}(\theta,\phi)$ be the far field pattern of $\Omega$ by the incident field $u^{inc}=e^{ik\db\cdot \bx }$. For a shifted domain $\Omega_h:=\{\bx+h,\bx\in \Omega\}$ with a constant vector $h\in \mathbb{R}^3$, the far field $u_{h}^{\infty}(\theta,\phi)$ generated by the incident wave $u^{inc}$ becomes
	\begin{eqnarray*}
	u_{h}^{\infty}(\theta,\phi) = e^{ik(\db-\hat{\bx})\cdot h}u^{\infty}(\theta,\phi)
	\end{eqnarray*}  
	where $\hat{\bx}=(\cos\theta\sin\phi, \sin\theta\sin\phi, \cos\phi)$.
\end{theorem}

\begin{proof}
	By shifting $\Omega_h$ back to $\Omega$ and making use of the uniqueness theorem of the exterior scattering problem \cite{Colton}, we see that the scattered field $u^\emph{scat}_{h}(\bx)$ and the normal derivative of the scattered field $\partial u^\emph{scat}_{h}(\bx)/\partial \nu(\bx)$ on $\partial \Omega_h $ are simply given by
	\begin{eqnarray*}
	\begin{cases}
	u^\emph{scat}_{h}(\bx) = u^\emph{scat}(\bx-h)e^{ik\db\cdot h}, \\
	\frac{\partial u^\emph{scat}_{h}(\bx)}{\partial \nu(\bx)} = \frac{\partial u^\emph{scat}(\bx-h)}{\partial \nu(\bx-h)}e^{ik\db\cdot h},
		\end{cases}
	\bx\in\Omega_h.		
	\end{eqnarray*}
	where $u^\emph{scat}(\bx-h)$ and $ \frac{\partial u^\emph{scat}(\bx-h)}{\partial \nu(\bx-h)}$ are respectively the scattered field and its normal derivative on $\Omega$. The conclusion now follows from equation \eqref{far_field} for the far field pattern.
	\end{proof}

Theorem \ref{thm:location} implies that the shifted domain of $\Omega$ simply changes the phase of the far field but not the modulus. Therefore, for a given far field data $u^{\infty}(\theta,\phi)$ of an unknown axis-symmetric obstacle $\Omega$, $|u^{\infty}(\theta,\phi)|$ is the same as the modulus of the far field of $\Omega_0$, where $\Omega_0$ denotes a shifted $\Omega$ with its axis centered at the origin. Thus by Theorem \ref{thm:rotation}, up to a rotation, $|u^{\infty}(\theta,\phi)|$ is symmetric with respect to $\theta=\theta_0$ and $\theta=2\pi-\theta_0$ for a given $\theta_0$ and any $\phi\in[0,\pi]$. This rotation angle is exactly the orientation of the axis of symmetry of the unknown obstacle $\Omega$. Once the axis is parallel to the $z$-axis, we can make use of the phase information to determine the $x$ and $y$ coordinates of the axis. In particular, by Theorem \ref{thm:location}, if we multiply the far field $u^{\infty}(\theta,\phi)$ by an appropriate factor $e^{ik\hat{x}\cdot h}$, the real and imaginary parts of the new far field will be symmetric with respect to  $\theta=\theta_0$ and $\theta=2\pi-\theta_0$ for a given $\theta_0$ and any $\phi\in[0,\pi]$. 

To summarize, we propose a three-step method to locate the axis of symmetry of the obstacle. In the first step, we evaluate the far field pattern based on the measured scattered field ${\bf u}^\emph{meas}_{k,\db}$ on the sphere $\partial\mathcal{B}$. Next, we determine the orientation of the axis of symmetry. In the third step, we find the location of its center. A detailed description of the algorithm follows:
\begin{enumerate}
	\item {\bf Step 1 (Evaluate the \textit{far field})}: For a fixed wavenumber $k$ and $\db$, collect the measured scattered field ${\bf u}^\emph{meas}_{k,\db}(\bx)$ on $\partial \mathcal{B}$ due to the incident plane wave $u^\emph{inc}(\bx)$. From the measurements of the scattered field on $\partial \mathcal{B}$, one can obtain the scattered field $u^\emph{scat}(\bx)$ anywhere on $\partial \mathcal{B}$ by using interpolation. To find the corresponding far field $u^{\infty}(\theta,\phi)$, we solve the boundary integral equation
	\begin{eqnarray*}
	\left(\frac{1}{2}\mathcal{I}+\mathcal{D}_{\partial \mathcal{B}}+ik\mathcal{S}_{\partial \mathcal{B}}\right)\eta(\bx)= u^{scat}_{k,\db}(\bx),\mbox{ on } \partial \mathcal{B}.
	\end{eqnarray*}
	where $\mathcal{D}_{\partial \mathcal{B}}$ and $\mathcal{S}_{\partial \mathcal{B}}$ are the single and double layer potentials defined on $\partial \mathcal{B}$. Once $\eta(\bx)$ is found, the far field can be evaluated using the formula
	\begin{eqnarray*}
	u^{\infty}(\theta,\phi) = \frac{1}{4\pi}\int_{\partial \mathcal{B}}\bigg\{\frac{\partial e^{-\mathrm{i}k \hat{\bx}\cdot \by}}{\partial \nu(\by)}+ike^{-\mathrm{i}k \hat{\bx}\cdot \by}\bigg\}\eta(\by)\mathrm{d}s(\by).
	\end{eqnarray*}
	\item {\bf Step 2 (Determine the \textit{orientation})}:  Suppose the far field $u^{\infty}(\theta,\phi)$ is given at $0=\theta_0< \theta_1<\cdots<\theta_i<\cdots<\theta_n=2\pi$ for $\theta$ and $0=\phi_0< \phi_1<\cdots<\phi_j<\cdots<\phi_m=\pi$ for $\phi$. Check if $|u^{\infty}(\theta,\phi)|$ is symmetric on the horizontal cross section by taking each $(\theta_i,\phi_j)$ as the north pole. If that is found within a certain accuracy, we take $(\theta_i,\phi_j)$ as the orientation of the axis of symmetry.
	\item {\bf Step 3 (Determine the \textit{location})}: Assume the orientation of the axis of symmetry is parallel to the $z$-axis.  In order to find the $x$ and $y$ coordinates of the axis, we send two incident waves $u_1^\emph{inc}(\bx)=e^{ikx}$ and $u_2^\emph{inc}(\bx)=e^{iky}$,  and evaluate their far field data $u_1^{\infty}(\theta,\phi)$ and $u_2^{\infty}(\theta,\phi)$ respectively. Suppose the obstacle is located in the area $[x_{\min},x_{\max}]\times [y_{\min},y_{\max}]$. By Theorem \ref{thm:location}, we can determine $h_1$ such that the real and imaginary parts of the shifted far field $u_2^{\infty}(\theta,\phi)e^{ik\cos\theta \sin \phi h_1}$ are symmetric with respect to $\theta=\pi/2$ and $\theta=3\pi/2$ for any $\phi\in[0,\pi]$ by doing an exhaustive search in the interval $[x_{\min},x_{\max}]$. Similarly, there exists  $h_2\in [y_{\min},y_{\max}]$ such that the real and imaginary parts of the shifted far field $u_1^{\infty}(\theta,\phi)e^{ik\sin\theta \sin \phi h_2}$ are symmetric with respect to $\theta=0$ and $\theta=\pi$ for any $\phi\in[0,\pi]$. In the end, we take $(h_1,h_2)$ as the $x$ and $y$ coordinates of the center of the obstacle.	
\end{enumerate}


\subsection{Recovering the shape of the obstacle(Part 2)} Once the axis of symmetry of the obstacle is obtained, we can reconstruct the shape of the obstacle. Our goal is to find an approximation of the simple open curve $\gamma:\left[0,1\right]\rightarrow\mathbb{R}^2$ that generates the surface of the obstacle using the measurements ${\bf u}_{k_j}^\emph{meas}$ for $j=1,\ldots,N_k$. The idea is to apply the RLA to solve a sequence of single frequency inverse scattering problems.

\subsubsection{Inverse scattering problem for a single frequency data} Using single frequency data, we can recast the inverse problem as the optimization problem
\begin{equation}
\tilde{\gamma} = \argmin_{\gamma} \|{\bf u}_k^\emph{meas}-\mathcal{F}_k(\gamma)\|, \label{eq:inv_prob_obs_3d}
\end{equation}
where $\tilde{\gamma}$ is an approximation of the curve $\gamma$.

The problem \eqref{eq:inv_prob_obs_3d} is both nonlinear and ill-posed. To treat the nonlinearity, we apply the iterative damped Gauss-Newton method. First, an initial guess for the approximation of the curve $\gamma$ is chosen, say $\gamma_0$. Next, in each step of this method, given an approximation $\gamma_j$ of the generating curve at the $j^{th}$ step, we update the curve to obtain $\gamma_{j+1}=\gamma_j+\alpha\delta \gamma$, with $\alpha>0$ being a chosen constant. To obtain $\delta \gamma$, we solve  
\begin{equation} \label{eq:newton_step}
\partial_{\gamma}\mathcal{F}_k(\gamma_j) \delta\gamma = {\bf u}_k^\emph{meas} -\mathcal{F}_k(\gamma_j),
\end{equation}
where $\partial_{\gamma} \mathcal{F}_k(\gamma_j)=\left[\partial_{\gamma} \mathcal{F}_{k,\db_1}(\gamma_j);\cdots;\partial_{\gamma} \mathcal{F}_{k,\db_{N_d}}(\gamma_j)\right]$ and $\partial_{\gamma} \mathcal{F}_{k,\db_i}(\gamma_j)$, $i=1,\ldots,N_d$ are, respectively, the Fr\'echet derivatives of $\mathcal{F}_k$ and $\mathcal{F}_{k,\db_i}$ with respect to $\gamma$ evaluated at the curve $\gamma_j$. The value of $v(\xb)=\partial_{\gamma}\mathcal{F}_{k,\db_i}(\gamma_j)\delta \gamma(\xb)$ is obtained by solving the Helmholtz equation 
\begin{equation*}
	\begin{cases}
\Delta v(\xb)+k^2 v(\xb) = 0 \quad \text{for} \quad \xb\in \mathbb{R}^3\setminus\overline{\Omega}_j,  \\
 v(\xb)=-k^2(\delta\gamma(\xb) \cdot \nu(\xb))\left[\frac{\partial u_{k,\db_i}}{\partial \nu}\right](\xb) \quad \text{for}\quad \xb\in\partial \Omega_j,
\end{cases}
\end{equation*}
$v$ satisfies the Sommerfeld radiation condition, $\Omega_j$ is the obstacle obtained by rotating the curve $\gamma_j$ around the axis of symmetry, $\nu(\xb)$ is the normal vector to the surface of $\Omega_j$ at $\xb$ and $\frac{\partial u_{k,\db_i}}{\partial \nu}$ is the normal derivative of the total field that is a solution for the problem \eqref{eq:helm} for the obstacle $\Omega_j$ with the incoming plane wave $u_{k,\db_i}^{inc}=e^{ik\xb\cdot\db_i}$. 

The iterations are repeated until a stopping criteria is reached. The stopping criteria can be the total number of iterations $N_{it}$, the residual achieving $\|{\bf u}_k^\emph{meas}-F_k(\gamma_j)\|\leq\epsilon_r$, with $\epsilon_r>0$, the difference of the curve evaluate in a set of points between consecutive steps is smaller than a certain value $\epsilon_s>0$, or others. A summary of the damped Gauss-Newton method is presented in Algorithm \eqref{alg:DGN}.

As mentioned, Problem \eqref{eq:inv_prob_obs_3d} is highly ill-posed. Various ways were proposed to deal with the ill-posedness of the inverse scattering problem, including Tykhonov regularization, truncated SVD \cite{Colton}, the use of a bandlimited representation of the domain \cite{Borges2015}, etc. In this work, we choose to search for a bandlimited representation of the generating curve $\gamma$. We represent the generating curve $\gamma$ and the update $\delta \gamma$ as 
\begin{equation}
\gamma(t)=p(t)\left(\cos\left(\pi \left(t-0.5\right)\right),\sin\left(\pi \left(t-0.5\right)\right)\right) \label{eq:rep_curve}
\end{equation}
and
\begin{equation}
\delta \gamma(t)=h(t)\left(\cos\left(\pi \left(t-0.5\right)\right),\sin\left(\pi \left(t-0.5\right)\right)\right) \label{eq:rep_curve1}
\end{equation}
where $p(t)$ and $h(t)$ are given by
\begin{equation}
p(t)=p^c_0+\sum_{j=1}^{N_p}\left(p^c_j\cos\left(2\pi j \left(t-0.5\right)\right)+p^s_j\sin\left(2\pi j \left(t-0.5\right)\right)\right)
\end{equation}
and
\begin{equation}
h(t)=h^c_0+\sum_{j=1}^{N_p}\left(h^c_j\cos\left(2\pi j \left(t-0.5\right)\right)+h^s_j\sin\left(2\pi j \left(t-0.5\right)\right)\right) \label{eq:rep_pol}
\end{equation}
with $p^c_j$, $p^s_j$, $h^c_j$ and $h^s_j$ being constant coefficients for the $j^{th}$ cosine and sine modes, $j=1,\ldots,N_p$. From Heisenberg's uncertainty principle for waves, we have that sub-wavelength features of the scatterer are present in the evanescent modes of the signal and are not detectable in finite precision. Consequently, the main advantage of this representation is that if we choose the bandlimit parameter $N_p=\mathcal{O}(k)$, the system of equations on \eqref{eq:newton_step} becomes well-conditioned. The second advantage of choosing this representation comes from the easy and fast evaluation of the polynomials $p(t)$ and $h(t)$ by using non-uniform FFT \cite{greengard2004accelerating, lee2005type}. The main disadvantage of choosing this representation stems from the fact that this representation is ideal for star-shaped figures. If the obstacle that we are trying to recover is not star-shaped, this representation will probably not work in terms of providing a high resolution reconstruction. An alternative is to use a bandlimited curve smoother like the one presented in \cite{Borges2015, beylkin2014fitting}.


\begin{algorithm}
\caption{Damped Gauss-Newton method}
\label{alg:DGN}
\begin{algorithmic}[1]
\STATE{{\bf Input:} Scattered field measurements ${\ub}^\emph{meas}_k$, initial guess $\gamma_0$, parameters $\alpha$, $N_{it}$, $\epsilon_r$, and $\epsilon_s$.}
\STATE{Set $j=0$, $\gamma=\gamma_0$ and $\gamma_{-1}=\gamma_0(1+2\|\gamma_0\|\epsilon_s)$.}
\WHILE{$j<N_{it}$ {\bf and}  $\|{\bf u}^\emph{meas}_k- \mathcal{F}_k(\gamma)\|\leq\epsilon_r$ {\bf and} $\|\gamma_j-\gamma_{j-1}\|\leq \epsilon_s$}
\STATE{Calculate $\mathcal{F}_k(\gamma_j)$ and $\partial_\gamma \mathcal{F}_k(\gamma_j)$.}
\STATE{Solve $\partial_{\gamma}\mathcal{F}_k(\gamma_j) \delta\gamma = {\bf u}_k^\emph{meas} -\mathcal{F}_k(\gamma_j)$.}
\STATE{$\gamma_{j+1}\leftarrow\gamma_j+\alpha \delta\gamma$}
\STATE{$j\leftarrow j+1$}
\ENDWHILE
\end{algorithmic}
\end{algorithm}

\subsubsection{Inverse scattering problem using multiple frequency data}
On the one hand, due to Heisenberg's uncertainty principle, the amount of information about the shape of the scatterer that can be stably recovered from measurements of the scattered field at frequency $k$ is proportional to $\mathcal{O}(k)$. This means that smaller features of the obstacle that have magnitude proportional to the sub-wavelength spectrum are extremely difficult to recover using finite precision. On the other hand, there are also inherent limitations to Newton-type methods for inverse scattering at a single frequency. When the incident field has larger wavelength, a low-resolution approximation of the inhomogeneity can be obtained using a simple initial guess. For problems with incident field with smaller wavelength, the initial guess for the iterative method must be close to the solution for the method to converge. This interplay between obtaining a low resolution reconstruction for large wavelengths using a simple initial guess, and the need to have a very good initial guess when using small wavelengths, together with the natural limitation on the amount of information that can be obtained using single frequency data led to the proposal of the RLA \cite{bao2015inverse, chen1995recursive, Chen}. 

In the RLA, a sequence of increasingly complicated nonlinear optimization problems like \eqref{eq:inv_prob_obs_3d} at successively higher frequencies is solved using a continuation path in frequency. At a given frequency $k$, one uses the damped Gauss-Newton method to solve the inverse problem \eqref{eq:inv_prob_obs_3d} and obtain an approximate solution $\gamma^{(k)}$ to the curve $\gamma$. This solution is used as the initial guess for the damped Gauss-Newton method to solve the nonlienar optimization problem with scattered field data at frequency $k+\delta k$, where $\delta k>0$ is sufficiently small. A summary of the RLA is presented in Algorithm \eqref{alg:RLA}.


\begin{algorithm}
\caption{Recursive Linearization Algorithm with Damped Gauss-Newton method}
\label{alg:RLA}
\begin{algorithmic}[1]
\STATE{{\bf Input:} Scattered field measurements ${\bf u}^\emph{meas}_{k_j}$, for $j=1,\ldots,N_k$ and $k_1<k_2<\ldots<k_{N_k}$, initial guess $\gamma_0$, parameters $\alpha$, $N_{it}$, $\epsilon_r$, and $\epsilon_s$.}
\STATE{Set $\gamma^{(0)}=\gamma_0$.}
\FOR{$j=1,\ldots,N_k$}
\STATE{Apply the damped Gauss-Newton method with initial $\gamma^{(j-1)}$, parameters $\alpha$, $N_{it}$, $\epsilon_r$, and $\epsilon_s$ to obtain an approximation $\tilde{\gamma}$ of the curve as a solution.}
\STATE{$\gamma^{(k)}\leftarrow\tilde{\gamma}$}.
\ENDFOR
\end{algorithmic}
\end{algorithm}



\section{Numerical Experiments}\label{s:num_res}

To illustrate our framework, we present five numerical examples. In Example 1, we show the results of part 1 of our framework by recovering the axis of symmetry of an oblique ellipsoid. In Example 2, we investigate the interplay between the frequency and the local sets of convexity of the objective functional $f_k(\gamma)=\|\bf{u}_k^\emph{meas}-\mathcal{F}_k(\gamma)\|$ when the obstacle is a sphere. In Example 3, we recover the shape of an obstacle using different geometric configurations regarding the direction of the incident plane wave and the position of the receptors. In Example 4, we recover the shape of the obstacle with different number of modes representing the domain. Finally, in Example 5, we recover the shape of an obstacle with sharp corners using the RLA. In this example in particular, we can see the effect of Gibbs phenomenon due to the approximation of the corners by a limited number of modes. To mitigate this effect, we apply a low-pass Gaussian filter to the update of the domain in each step. In Examples 3, 4 and 5, we consider that we have already applied part 1 of our framework and that we have the axis of symmetry and location of the obstacle with high precision. Hence, we present only the results of part 2 of the framework in those examples.

A list of the numerical examples with their respective descriptions and results are presented in Table \ref{tab:examples}. 
\begin{table}[!htbp]
	\caption{List of numerical examples with respective tables and figures.}\label{tab:examples}
	\begin{center}
		\begin{tabular}{l  l cc}\hline
			Example &  Description & Tables & Figures \\ \hline
			1 & Recovering the axis of symmetry & X & \ref{fig:ex1_orientation}, \ref{fig:ex1_location} \\
			2 & Interplay between frequency and initial guess for a sphere & \ref{tab:example2} & \ref{fig:example2} \\
			3 & Illumination of the obstacles and placement of the receptors & X & \ref{fig:ex3_general}, \ref{fig:ex3_rec}\\
			4 & Limiting the number of modes used to recover the obstacle & X & \ref{fig:ex4_rec}\\
			5 & Reconstruction of sharp features using multiple frequencies & X & \ref{fig:ex5_general}, \ref{fig:ex5_rec}\\\hline
		\end{tabular}
	\end{center}
\end{table}

\subsection{Example 1: Recovering the axis of symmetry} \label{ex:1}
\begin{figure}
	\center
	\begin{subfigure}{.32\textwidth}
		\includegraphics[scale=0.24]{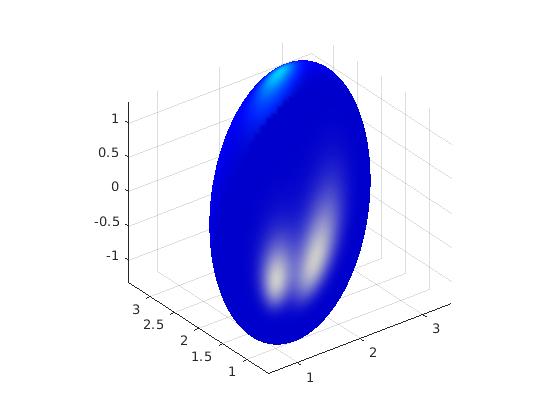}
		\caption{An oblique ellipsoid}\label{fig:rotated_ellip}
	\end{subfigure}
	\begin{subfigure}{.32\textwidth}
		\includegraphics[scale=0.24]{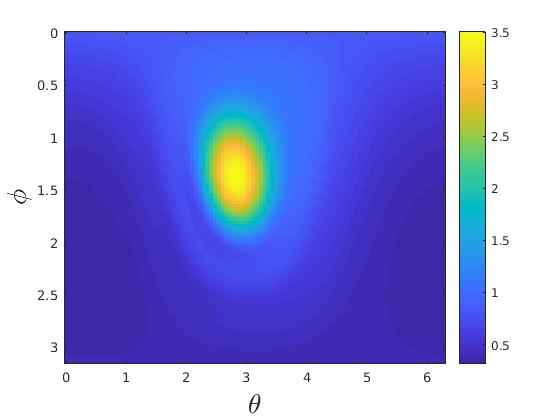}
		\caption{$\lvert u^{\infty}\rvert$ before rotation.}\label{fig:ex1_far1}
	\end{subfigure}
	\begin{subfigure}{.32\textwidth}
		\includegraphics[scale=0.24]{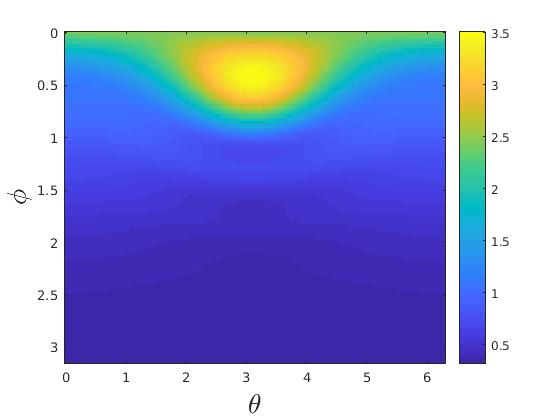}
		\caption{$\lvert u^{\infty}\rvert$ after rotation. }\label{fig:ex1_far2}
	\end{subfigure}
	\caption{{\bf (Example 1)} We present in (a) an oblique ellipsoid centered at $(2,2,0)$ , (b) the modulus of the far field pattern before rotation, and (c) the modulus of the far field pattern after rotation with angel given by $\theta = 0.25\pi$, $\phi =0.33\pi$. }\label{fig:ex1_orientation}
\end{figure}

\begin{figure}
	\center
	\begin{subfigure}{.32\textwidth}
		\includegraphics[scale=0.24]{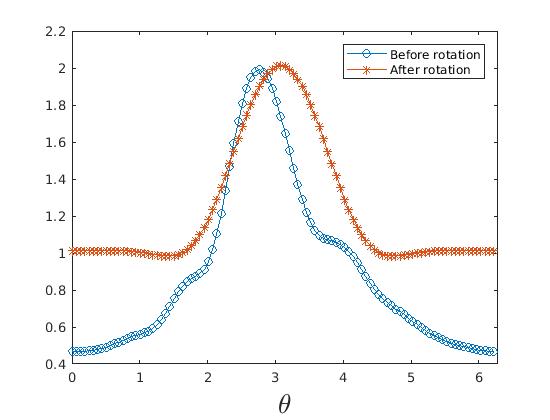}
		\caption{ Comparison of $\lvert u^{\infty}\rvert$ before and after rotation}\label{fig:ex1_ori}
	\end{subfigure}
	\begin{subfigure}{.32\textwidth}
		\includegraphics[scale=0.24]{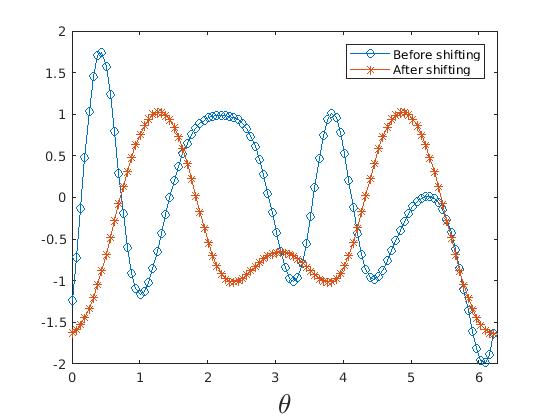}
		\caption{Comparison of $ u^{\infty}$ before and after shifting}\label{fig:ex1_loc_x}
	\end{subfigure}
	\begin{subfigure}{.32\textwidth}
		\includegraphics[scale=0.24]{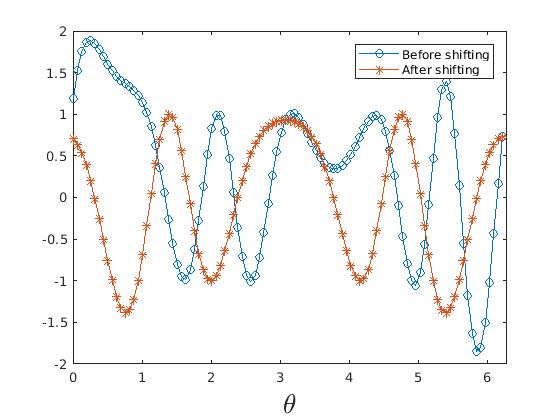}
		\caption{Comparison of $ u^{\infty}$ before and after shifting}\label{fig:ex1_loc_y}
	\end{subfigure}
	\caption{{\bf (Example 1)} We have in (a) the comparison of the modulus of the far field pattern at $\phi =\pi/4$ before and after rotation, in (b) the comparison of the real part of the far field pattern at $\phi =\pi/4$ before and after the shifting along the $x$ direction, and in (c) the comparison of the real part of the far field pattern at $\phi =\pi/4$ before and after the shifting along the $y$ direction.}\label{fig:ex1_location}
\end{figure}

We test part 1 of our framework by determining the orientation and location of the axis of symmetry of an unknown obstacle. In particular,  we are trying to determine the axis of an oblique ellipsoid generated by rotating and shifting from a standard one, whose parameterization of generating curve is given by 
\begin{eqnarray}
\gamma(t)=\left(x(t),y(t)\right)= \left(\cos(t-0.5),  2\sin(t-0.5)\right), t\in [0,\pi]. \nonumber
\end{eqnarray}
The axis of the ellipsoid is oriented at $\theta = \frac{\pi}{4}$, $\phi =\frac{\pi}{3}$  with the center located at $(2,2,0)$, as shown in Figure \ref{fig:rotated_ellip}. Using the incident wave  $u^\emph{inc}(\bx)=e^{ik\bx\cdot \db}$ with $k=3$ and  $\db=(\cos(\pi)\sin(\pi/8)$, $\sin(\pi)\sin(\pi/8),\cos(\pi/8))$, we measure the scattered field $u^\emph{scat}(\bx)$ on $\partial \mathcal{B}$ with $\mathcal{B}$ being radius $5$ and centered at the origin. The far field pattern $u^{\infty}(\hat{\bx})$ is found according to step 1 in part 1 of our algorithm. The modulus of $u^{\infty}(\hat{\bx})$ is shown in Figure \ref{fig:ex1_far1}. In particular, the symmetry can not be seen directly. We therefore test different rotation angles as the north pole as stated in the step 2 of part 1 by choosing $m=n=100$. Applying this test, we successfully find the orientation angle $\theta = 0.25\pi$ and $\phi =0.33\pi$, which is very close to the exact solution. The far field pattern after rotation is shown in Figure \ref{fig:ex1_far2}. We also plot the cross section of the far field pattern at $\phi= \pi/4$ before and after the rotation in 
Figure \ref{fig:ex1_ori}. One clearly sees that the symmetry is recovered if the correct rotation is found.  

Once the orientation is found, our next step is to determine the location of the axis on the $xy-$plane. We collect the far field data by sending two incident waves, respectively. One is $e^{ikx}$ and the other is $e^{iky}$. Since the location is not at the center, the far field patterns are not symmetric anymore. However, by determining the corresponding phase function $e^{i\alpha}$ and $e^{i\beta}$ from step 3 of part 1, where $\alpha=kh_1\cos\theta \sin \phi$, $\beta=kh_2\sin\theta\sin \phi$, we are able to see the symmetry of the far field pattern, as shown in Figures \ref{fig:ex1_loc_x} and \ref{fig:ex1_loc_y} for a cross section at $\phi=\pi/4$. Based on this fact, we recover the location of the center on the $xy-$plane as $(2,2)$.

\subsection{Example 2: Interplay between frequency and initial guess for a sphere} \label{ex:2}

In this example, let the objective functional be $f_k(\gamma)=\|{\bf u}_k^\emph{meas}-\mathcal{F}_k(\gamma)\|$. We consider the case where the obstacle is a sphere of radius $1$, with generating curve $\gamma:[0,1]\rightarrow\mathbb{R}^2$ given by $\gamma(t)=(\cos(\pi(t-0.5),\sin(\pi(t-0.5)))$. Applying the damped Gauss-Newton method one tries to recover at each step a polynomial $p(t)=p_0$. The shape reconstruction problem of finding $p_0$ turns into the problem of finding the root of a single variable nonlinear equation. To guarantee the convergence of the damped Gauss-Newton method, the initial guess must be in the same local convexity set of $f_k$ as the solution. This example aims to illustrate the interplay between the wavenumber of the incident wave and the local set of convexity near the solution. 

We use for each experiment one incident plane wave with incident direction $\db=(\cos(\pi/9)$, $0$, $\sin(\pi/9))$ and wavenumber $k_m$, such that $k_1=1$, $k_m=2.5(m-1)$, $m=2,\ldots,13$. The scattered data is measured at receptors $x_{\theta,\phi}=\rho\left(\cos(\phi_j) \sin(\theta_l),\sin(\phi_j)\sin(\theta_l),\cos(\theta_l)\right)$, with $\rho=10$, $\theta_l=l\pi/11$, $\phi_j=2j\pi/10$ for $l,j=1,\ldots,10$. Since we want to show the relation between the wavenumber and the local convexity set of $f_k$, we do not add noise to the measurements.

We calculate $f_{k_m}(\tilde{\gamma}^{(j)})$ at the curves $\tilde{\gamma}^{(j)}(t)=p^{(j)}_0(\cos(t),\sin(t))$, where $p^{(j)}_0=0.01j$, for $j=1,\ldots,1000$. Using the values of the objective functional, we are able to identify the local set of convexity in which $p_0=1$ is located. We denote this set to be $\left[a,b\right]$, where $0<a$ is the largest point smaller than $1$ where $f_k$ attains a local maximum, and $b$ is the smallest point larger than $1$ where the function attains a local maximum.

\begin{figure}[ht]
  \center
\begin{subfigure}[t]{0.48\textwidth}
\includegraphics[height=2.5in]{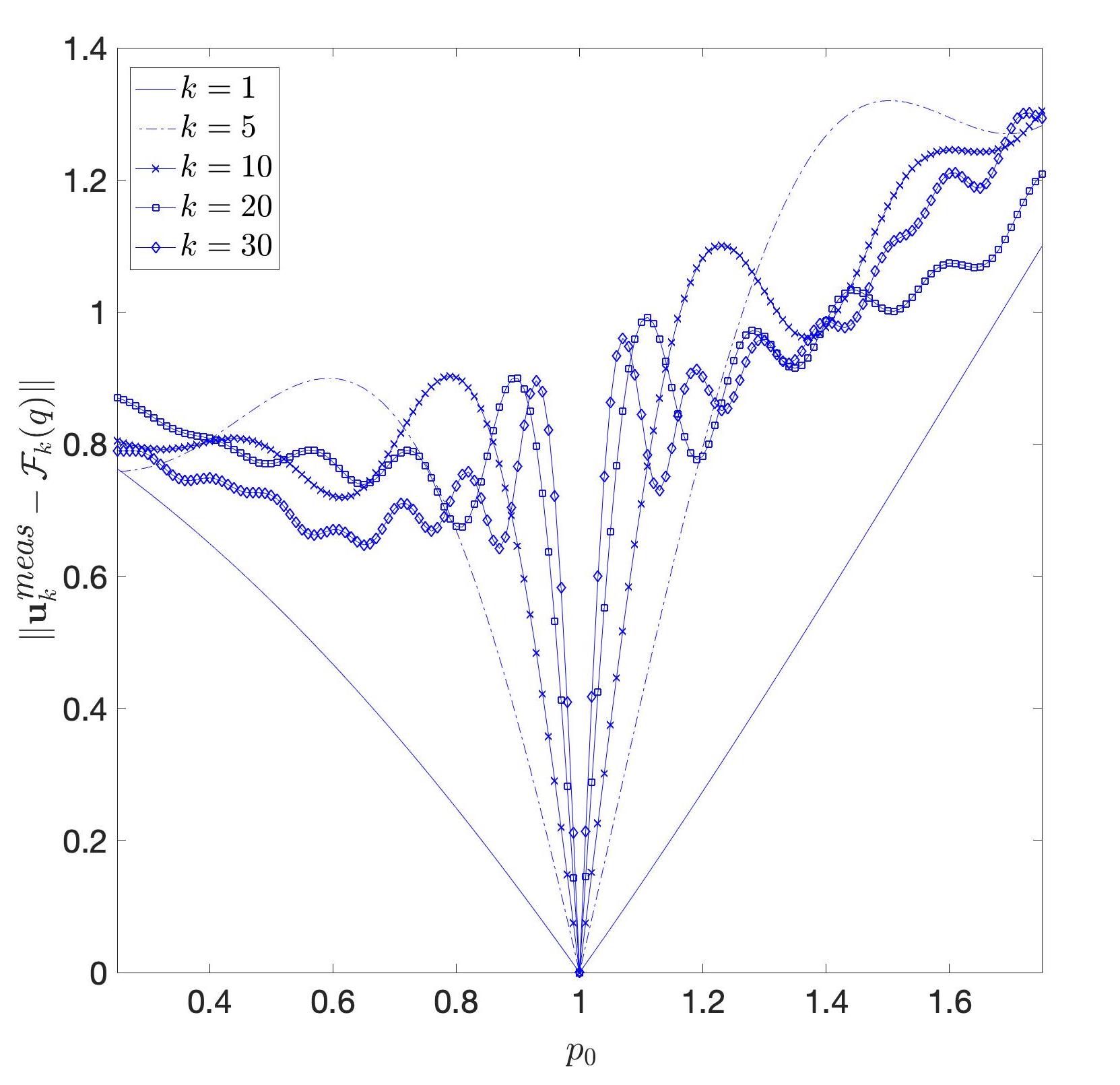} 
\caption{Objective functional $\|{\bf u}_k^\emph{meas}-\mathcal{F}_k(\gamma)\|$.}\label{fig:example2a}
\end{subfigure}
\begin{subfigure}[t]{0.48\textwidth}
\includegraphics[height=2.48in]{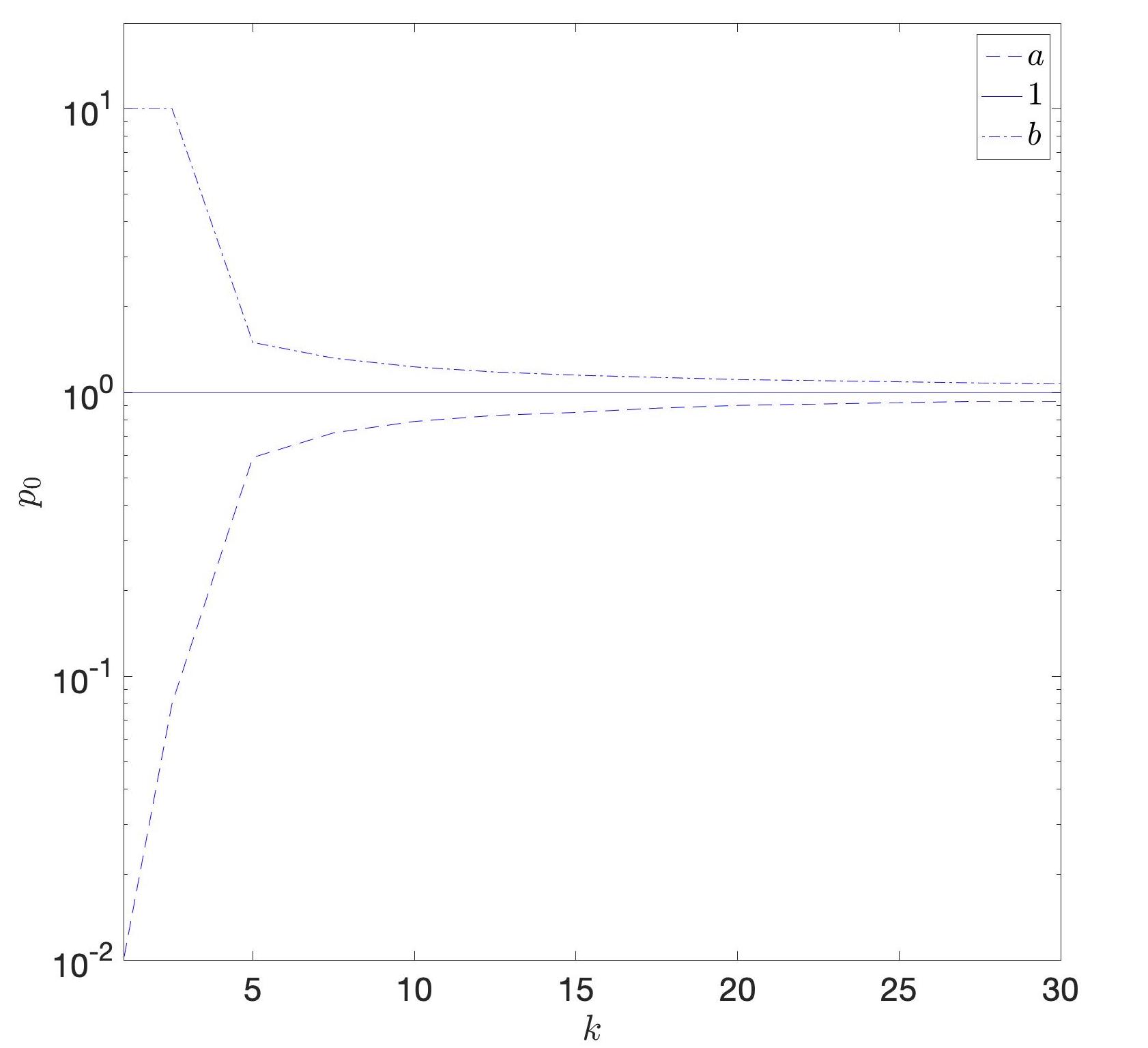} 
\caption{Local set of convexity $\left[a,b\right]$.}
\label{fig:example2b}
\end{subfigure}
  \caption{{\bf (Example 2)} We present in (a) the plot of the objective functional for different $k$, and in (b) the values of $a$ and $b$, which are the endpoints of the local set of convexity containing the root $p_0=1$.} \label{fig:example2} 
\end{figure}

In Figure \ref{fig:example2a}, we present the value of $f_{k}(\tilde{\gamma}^{(j)})$ for $k=1$, $5$, $10$, $20$ and $30$. In Figure \eqref{fig:example2b}, we plot three lines: the line for $p_0=1$, the line with values of $b$ and the line with values of $a$. The values of $a$ and $b$ are also available in Table \ref{tab:example2} with the wavenumber and its respective wavelength $\lambda=2\pi/k$.

\begin{table}
  \caption{Table with the values of the boundaries of the local set of convexity $\left[a,b\right]$ for different values of the wavenumber $k$ with the respective wavelength $\lambda=2\pi/k$.} \label{tab:example2} 
\begin{tabular}{c || *{13}{c}}
$k$            & 1.0   &  2.5    & 5.0   & 7.5   & 10.0   & 12.5   & 15.0   & 17.5   & 20.0   & 22.5   & 25.0   & 27.5   & 30.0 \\\hline
$\lambda$ & 6.28 &  2.51  & 1.26 & 0.84 &   0.63 &   0.50 &   0.42 &   0.36 &   0.31 &   0.28 &   0.25 &   0.23 &   0.21 \\ \hline\hline
$a$            & 0.01 &  0.08  & 0.59 & 0.72 &   0.79 &   0.83 &   0.85 &   0.88 &   0.90 &   0.91 &   0.92 &   0.93 &   0.93    \\
$b$            & 10.0 &  10.0  & 1.50 & 1.32 &   1.23 &   1.18 &   1.15 &   1.13 &   1.11 &   1.10 &   1.09 &    1.08 &   1.07     \\
\end{tabular}
\end{table}

As expected, with the increasing value of the wavenumber $k$, the size of the interval $\left[a,b\right]$ decreases. Also, as the wavenumber increases, $f_k$ presents multiple local minima, which shows the nonlinearity of the inverse problem. 

\subsection{Example 3: Illumination of the obstacle and placement of the receptors} \label{ex:3}
\begin{figure}
\center
\begin{subfigure}{.33\textwidth}
\includegraphics[scale=0.12]{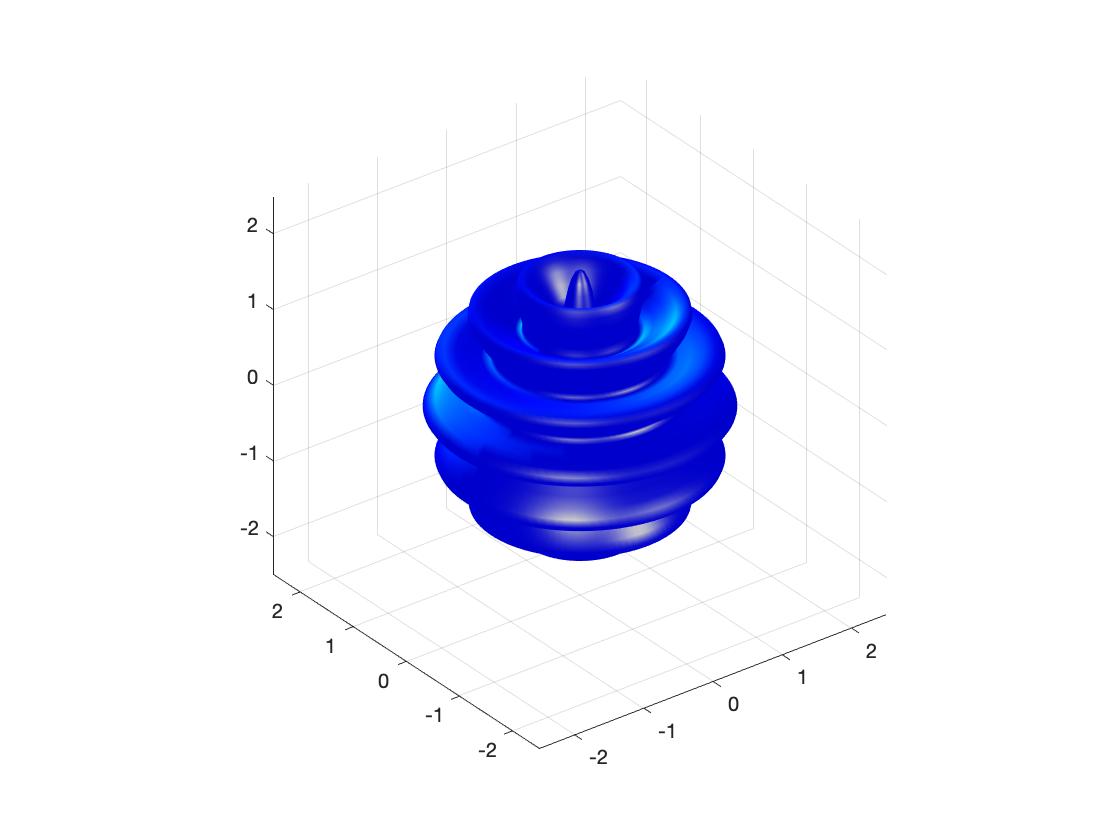}
\caption{Original obstacle}\label{fig:ex3_orig}
\end{subfigure}
\begin{subfigure}{.31\textwidth}
\includegraphics[scale=0.11]{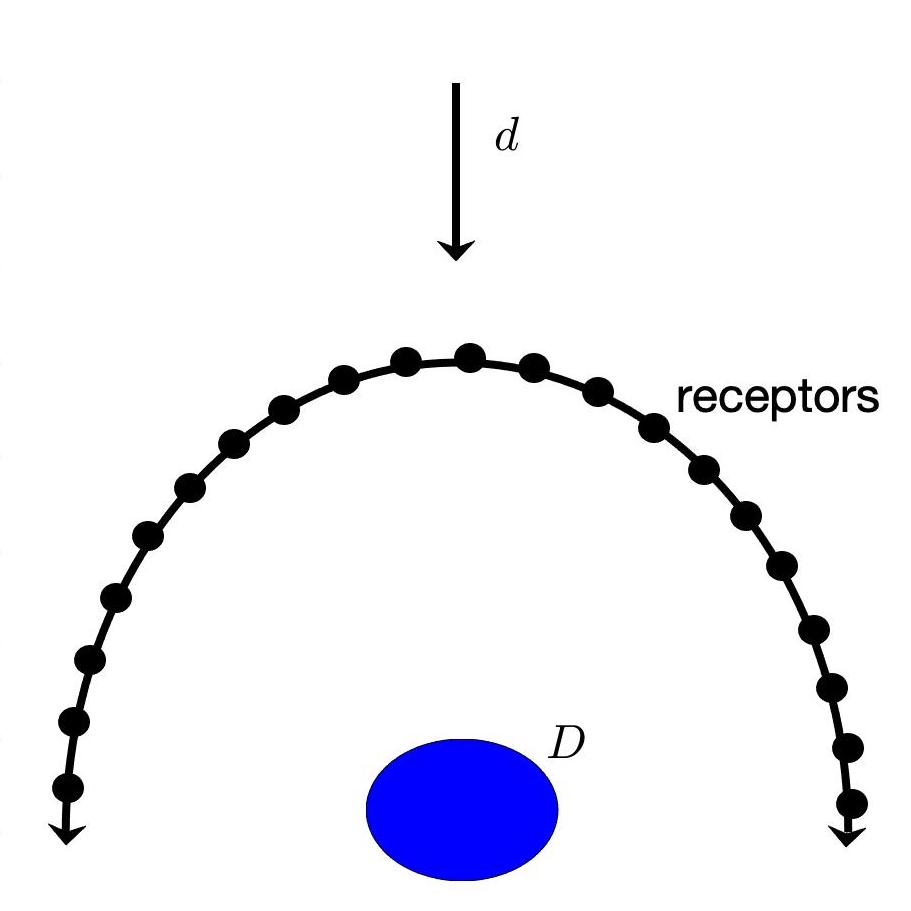}
\caption{Configurations 1 and 3}\label{fig:ex3_conf_1_3}
\end{subfigure}
\begin{subfigure}{.31\textwidth}
\includegraphics[scale=0.08]{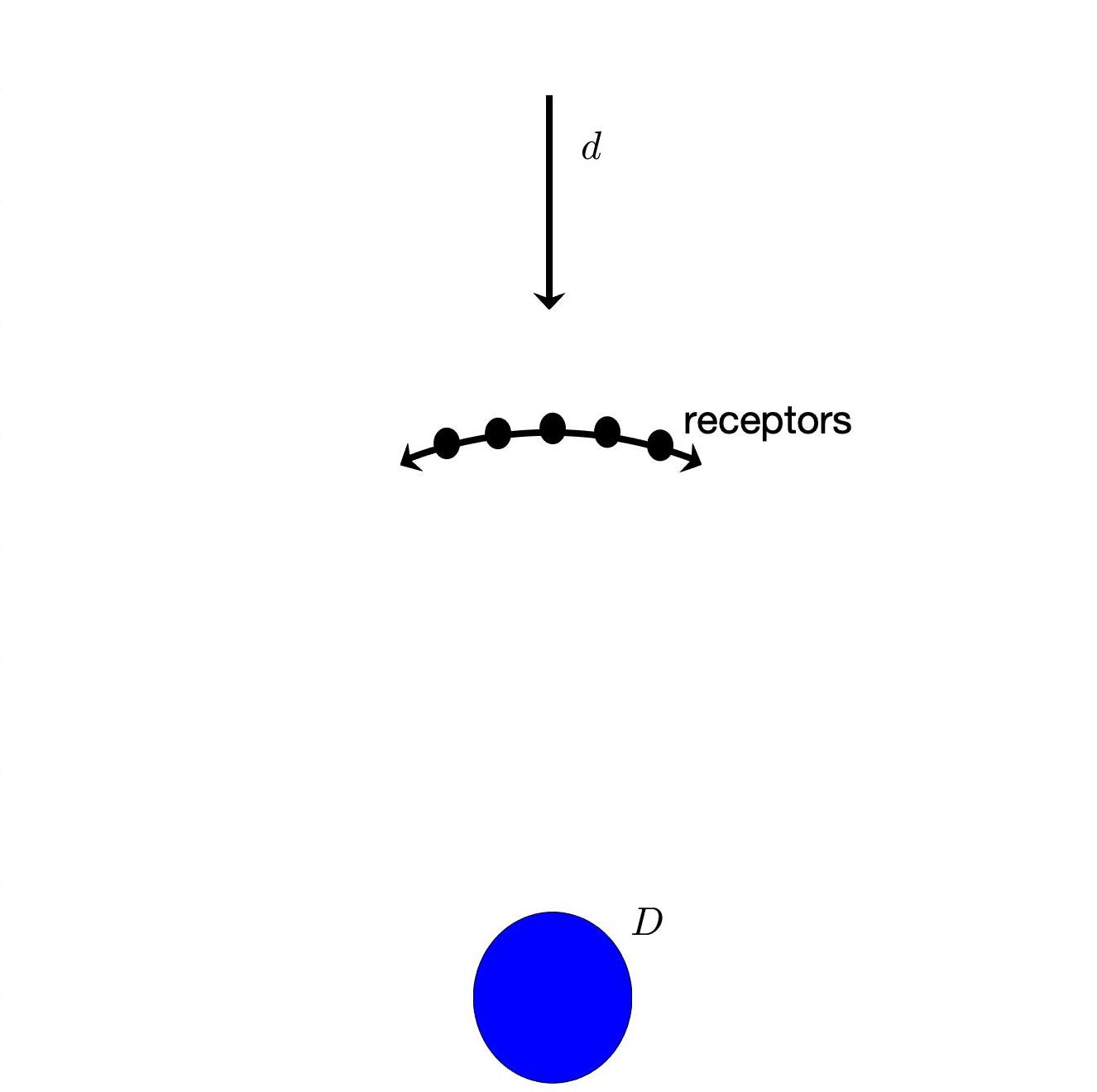}
\caption{Configurations 2 and 4}\label{fig:ex3_conf_2_4}
\end{subfigure}
\caption{{\bf (Example 3)} We present in (a) the original obstacle in Examples 3 and 4, (b) the sketch of configurations 1 and 3, and (c) the sketch of configurations 2 and 4.}\label{fig:ex3_general}
\end{figure}

In this example, we use multifrequency scattered data generated by four different geometric configurations of incident waves and receptors to recover the shape of an obstacle generated by the rotation around the $z$-axis of the curve $\gamma:[0,1]\rightarrow\mathbb{R}^2$, given by $\gamma(t)=p(t)\left(\cos\left(\pi\left(t-0.5\right)\right),\sin\left(\pi\left(t-0.5\right)\right)\right)$, where
\begin{equation*}
p(t)=1.5+(0.3\cos(8\pi(t-0.5)/0.5).
\end{equation*}
A three-dimensional rendering of the obstacle can be seen in Figure \ref{fig:ex3_orig}.

We set up the incident wave directions and the receptors positions in four different geometric configurations:
\begin{enumerate}
\item The incoming direction of the incident waves is $\db=(-1,0,0)$. The scattered field is measured at $N_r=100$ receptors located at the points $\xb_m=10\left(\cos(\theta_m),0,\sin(\theta_m)\right)$, with $\theta_m=-\pi/2+m\pi/(N_r-1)$, $m=0,\ldots,N_r-1$. See Figure \ref{fig:ex3_conf_1_3}.
\item The incoming direction of the incident waves is $\db=(-1,0,0)$. The scattered field is measured at $N_r=100$ receptors located at the points $\xb_m=10\left(\cos(\theta_m),0,\sin(\theta_m)\right)$, with $\theta_m=-\pi/2+((N_\theta-N_r)/2+m)\pi/(N_\theta-1)$, $m=0,\ldots,N_r-1$, and $N_\theta=3100$. See Figure \ref{fig:ex3_conf_2_4}.
\item The incoming direction of the incident waves is $\db=(0,0,-1)$. The scattered field is measured at $N_r=100$ receptors located at the points $\xb_m=10\left(\cos(\theta_m),0,\sin(\theta_m)\right)$, with $\theta_m=m\pi/(N_r-1)$, $m=0,\ldots,N_r-1$. See Figure \ref{fig:ex3_conf_1_3}.
\item The incoming direction of the incident waves is $\db=(0,0,-1)$. The scattered field is measured at $N_r=100$ receptors located at the points $\xb_m=10\left(\cos(\theta_m),0,\sin(\theta_m)\right)$, with $\theta_m=((N_\theta-N_r)/2+m)\pi/(N_\theta-1)$, $m=0,\ldots,N_r-1$ and $N_\theta=3100$. See Figure \ref{fig:ex3_conf_2_4}.
\end{enumerate}

In the configuration 1, the incident wave illuminates the obstacle in the direction perpendicular to the axis of symmetry and the receptors are located such that is possible to see the entire obstacle from their positions (taking symmetry into consideration). In configurations 2, 3 and 4, the position of the receptors provide only limited information about the obstacle.

To obtain measurements, we start by computing the scattered field data $u^\text{scat}_{k,\db}$ at frequencies $k_q=0.5+(q-1)0.25$, with $q=1,\ldots,25$. Next, to avoid inverse crimes, we add 2\% noise to the measured scattered data, using the formula
\begin{equation}
({\bf u}_{k,\db}^\emph{meas})_{m+1} = u^\emph{scat}(\xb_m) + 0.02 \frac{\epsilon}{\|\epsilon\|}|u_{k,\db}^\emph{scat}(\xb_m)| \label{eq:scat_noise}
\end{equation}
where $({\bf u}_{k,\db}^\emph{meas})_{m+1}$ is the $(m+1)^{th}$ coordinate of the vector ${\bf u}_{k,\db}^\emph{meas}$, $m=0,\ldots,N_r-1$, $\epsilon=\epsilon_1+i\epsilon_2$, and $\epsilon_1$ and $\epsilon_2$ are chosen from the random normal distribution $\mathcal{N}(0,1)$ with mean zero and variance one.

We apply the RLA with the damped Gauss-Newton method at each frequency. We set the stopping criteria of the Gauss-Newton method to be the maximum number of iterations $N_{it}=10$, the update step size should be no smaller than $\epsilon_s=0.03$ and the residual smaller than $\epsilon_r=0.03$. We also include as a stopping criteria any residual increase from one step to another. We set the damping parameter $\alpha=0.1$ for the damped Gauss-Newton method at frequency $k=0.5$, and $\alpha=0.1/\|h\|$ for all other frequencies, where $\|h\|$ is the 2-norm of the vector of coefficients of the update. We set the the number of modes in the polynomial representing the update $h$ to be $N_p=\lfloor2 k\rfloor$.  

\begin{figure}
\center
\begin{subfigure}{.45\textwidth}
\includegraphics[scale=0.2]{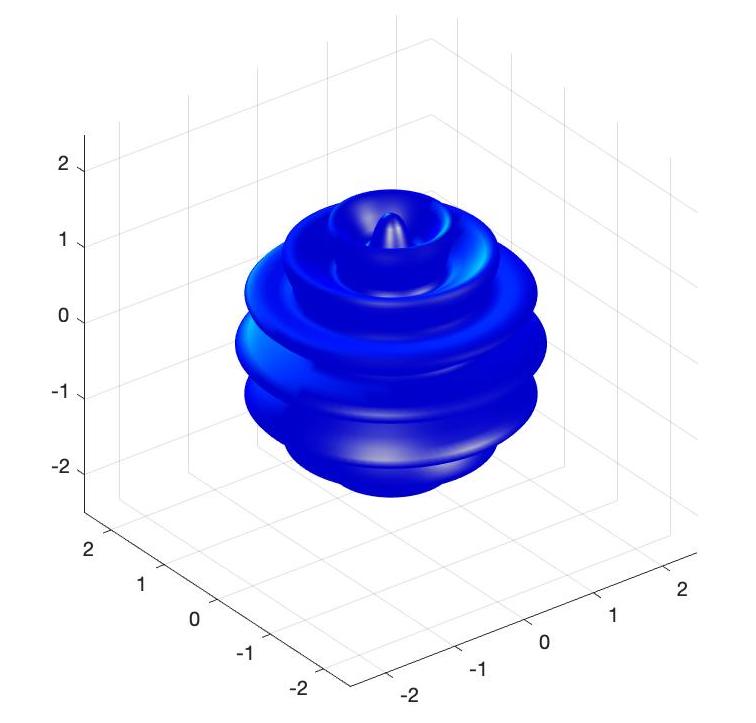}
\caption{Configuration 1 }\label{fig:ex3_conf1_3d_65}
\end{subfigure}
\begin{subfigure}{.45\textwidth}
\includegraphics[scale=0.2]{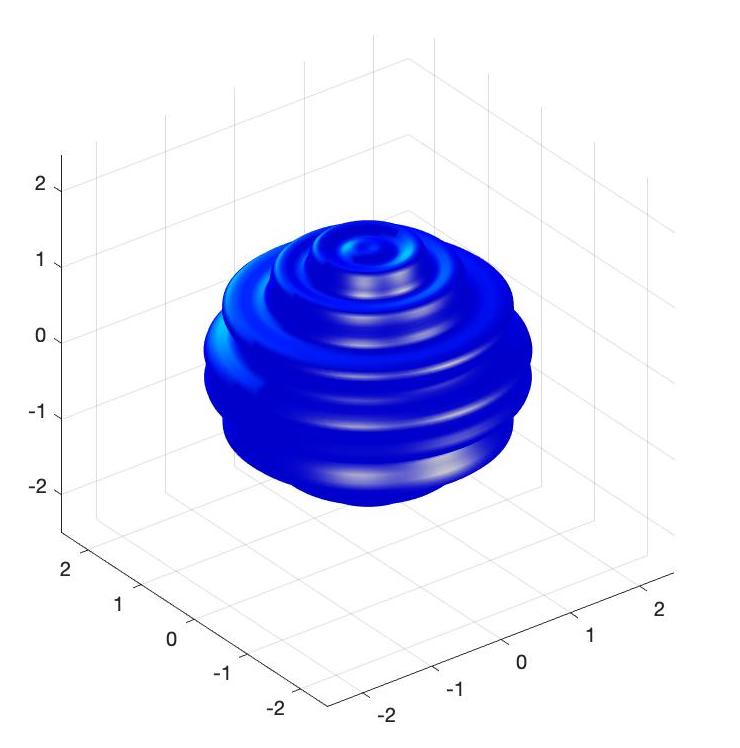}
\caption{Configuration 2}\label{fig:ex3_conf2_3d_65}
\end{subfigure}

\begin{subfigure}{.45\textwidth}
\includegraphics[scale=0.2]{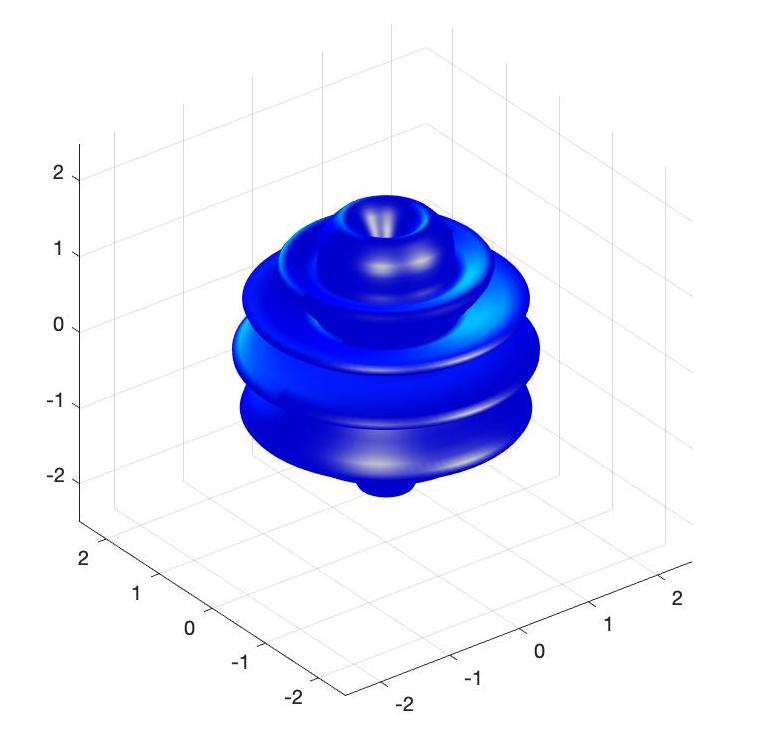}
\caption{Configuration 3}\label{fig:ex3_conf3_3d_65}
\end{subfigure}
\begin{subfigure}{.45\textwidth}
\includegraphics[scale=0.2]{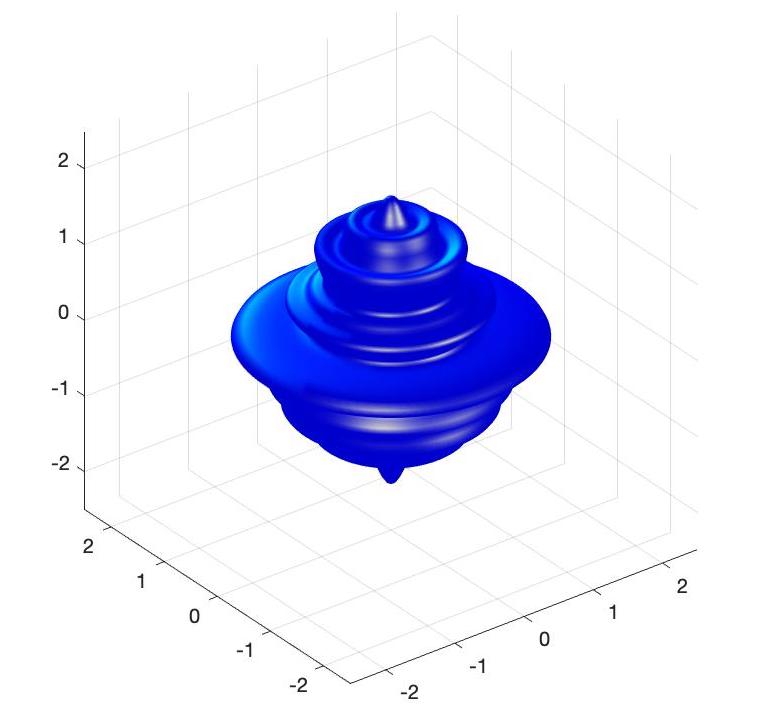}
\caption{Configuration 4}\label{fig:ex3_conf4_3d_65}
\end{subfigure}
\caption{{\bf (Example 3)} Reconstruction of the shape of the obstacle at wavenumber $k=6.5$ for different configurations.}\label{fig:ex3_rec}
\end{figure}

\begin{figure}
\center
\begin{subfigure}{.45\textwidth}
\includegraphics[width=.9\textwidth]{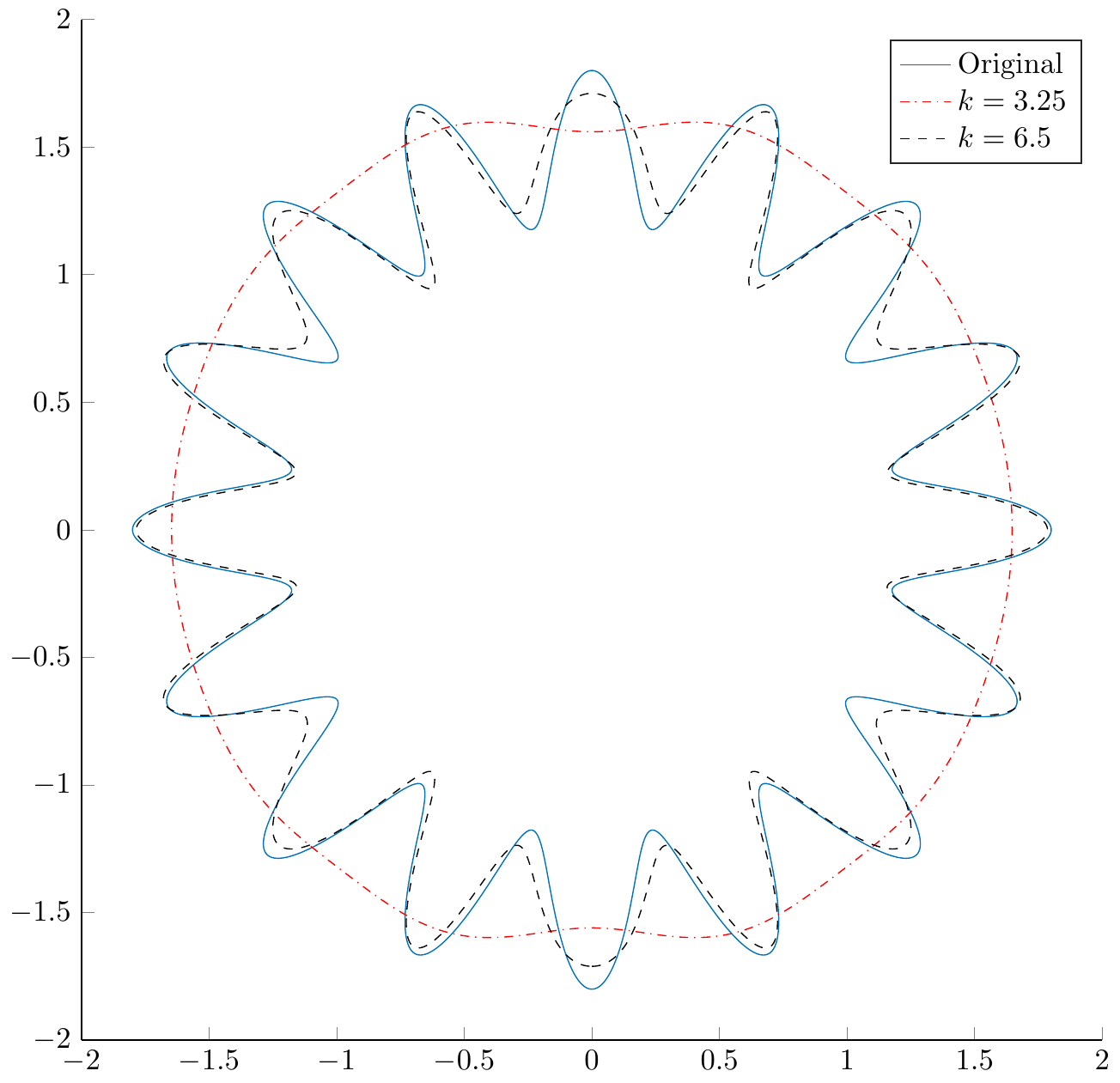}
\caption{Configuration 1}\label{fig:ex3_conf1_cs}
\end{subfigure}
\begin{subfigure}{.45\textwidth}
\includegraphics[width=.9\textwidth]{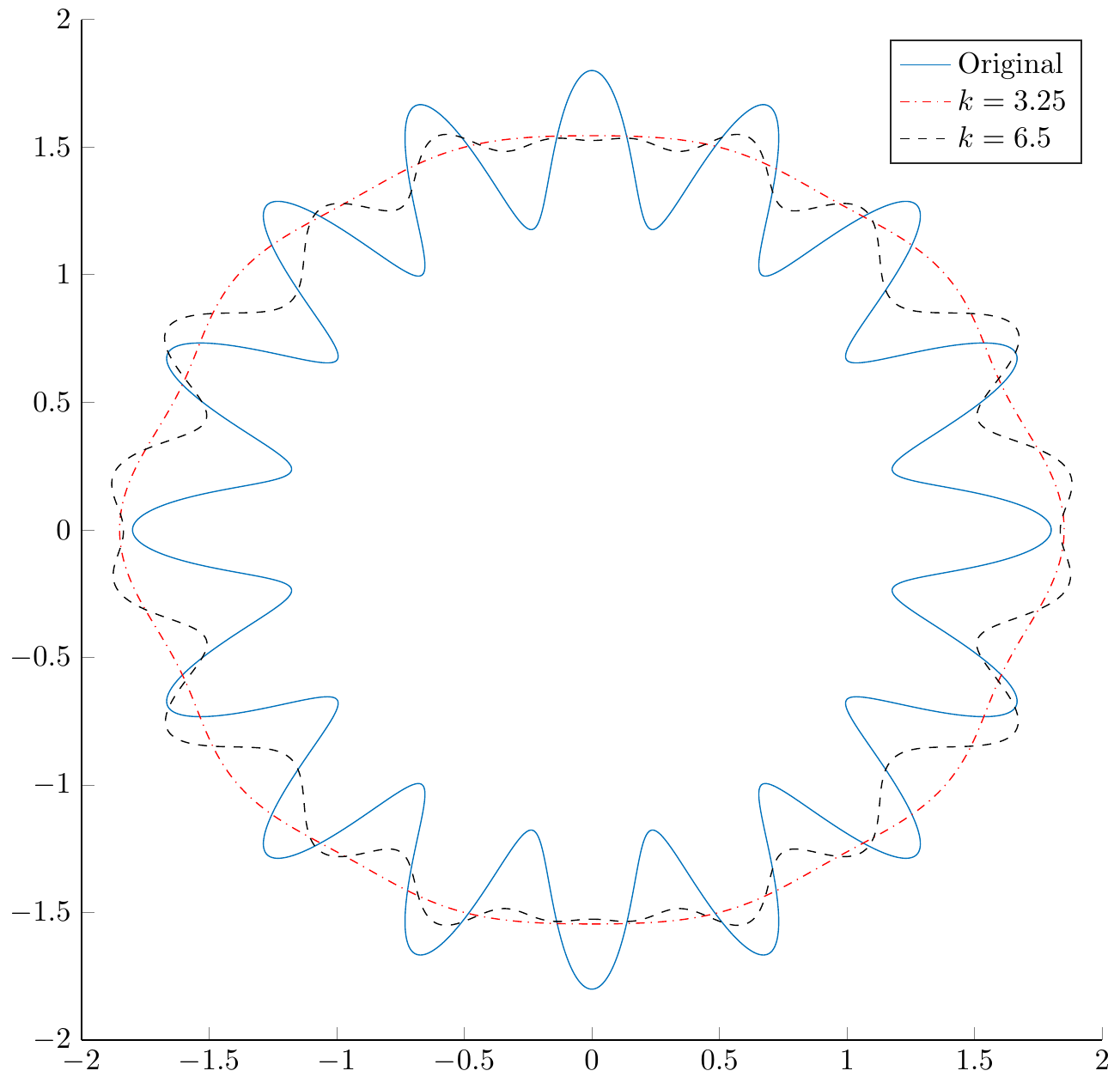}
\caption{Configuration 2}\label{fig:ex3_conf2_cs}
\end{subfigure}

\begin{subfigure}{.45\textwidth}
\includegraphics[width=.9\textwidth]{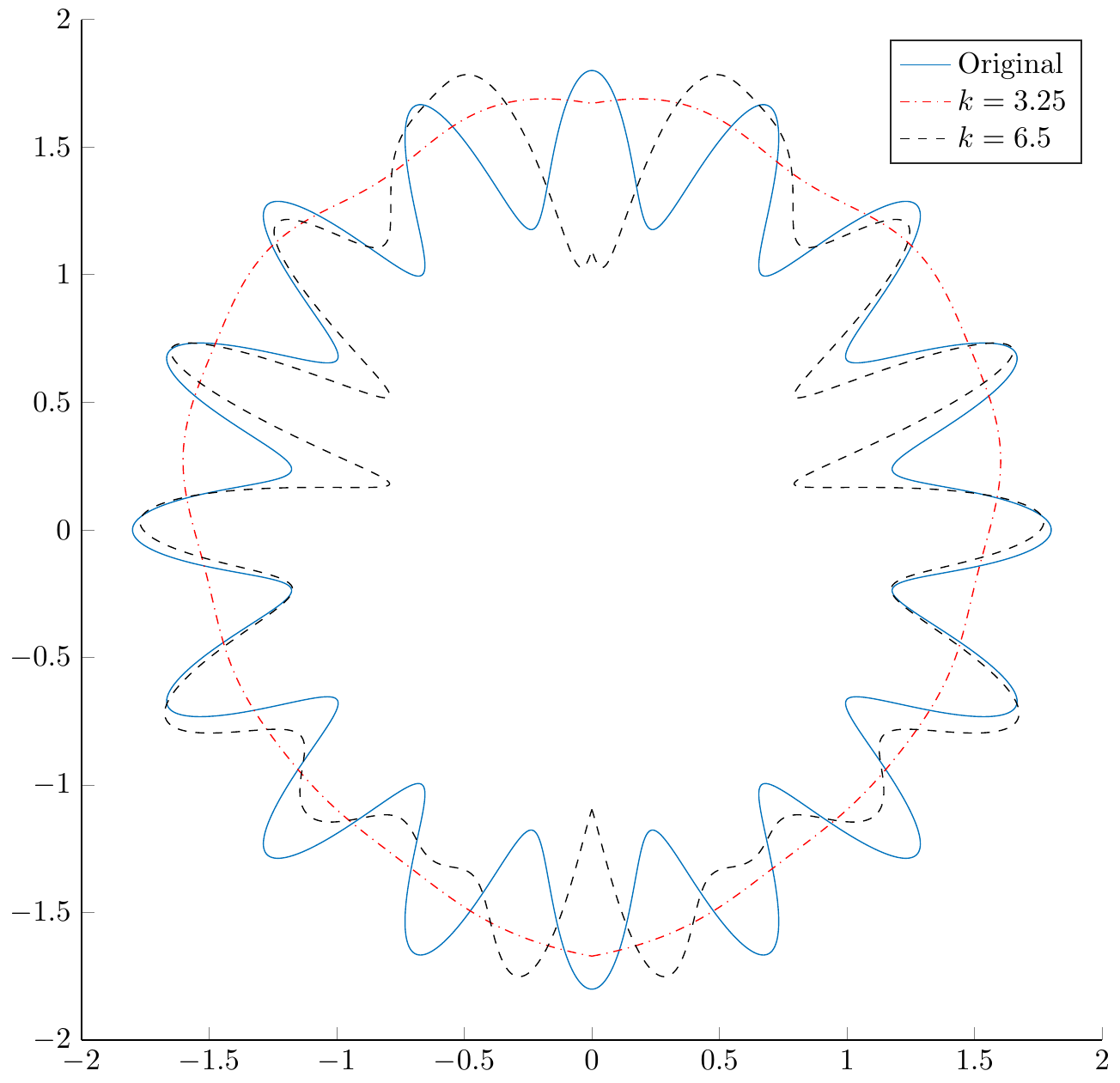}
\caption{Configuration 3}\label{fig:ex3_conf3_cs}
\end{subfigure}
\begin{subfigure}{.45\textwidth}
\includegraphics[width=.9\textwidth]{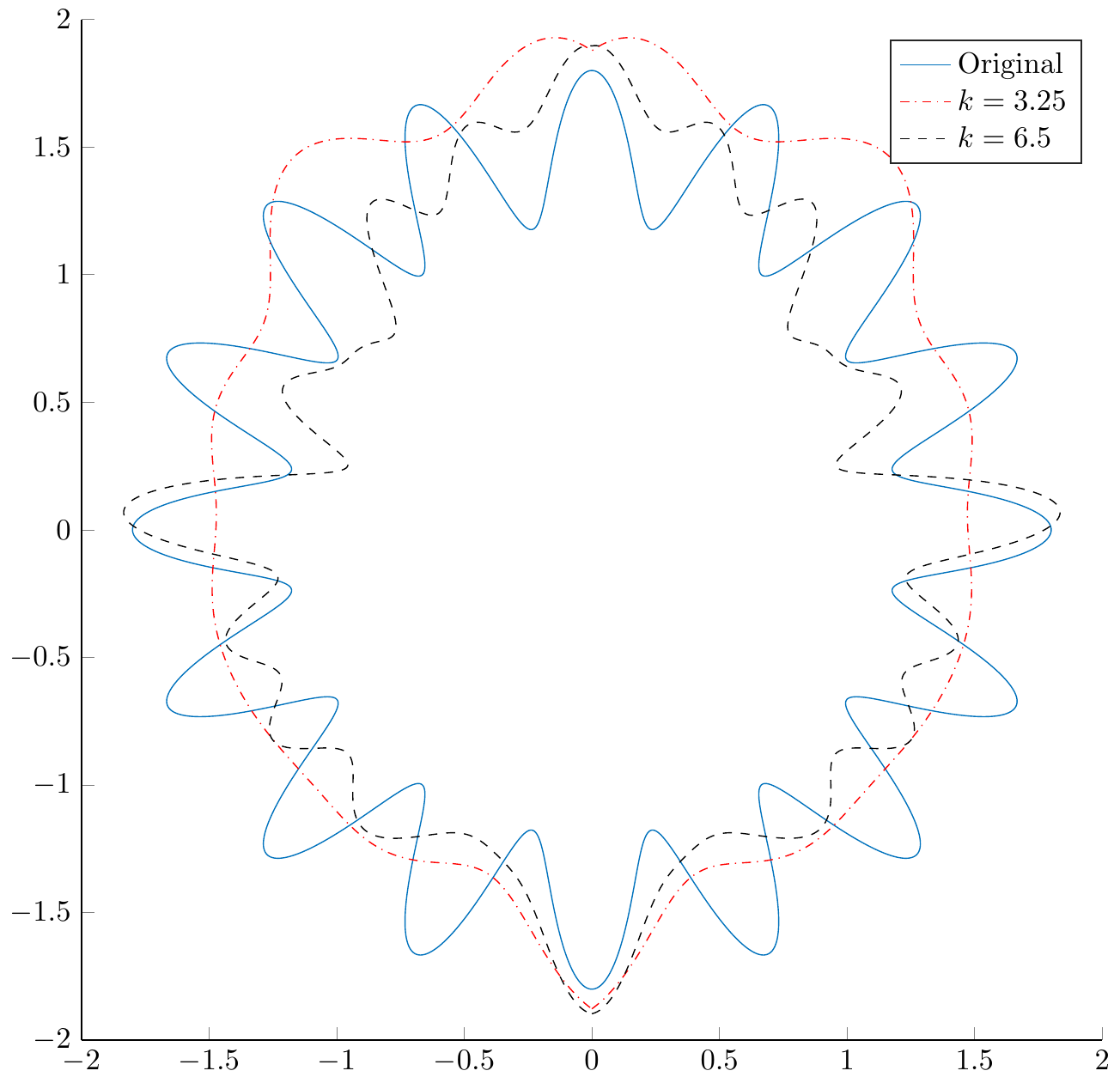}
\caption{Configuration 4}\label{fig:ex3_conf4_cs}
\end{subfigure}
\caption{{\bf (Example 3)} Cross section of the original obstacle with the reconstructions at $k=3.25$ and $6.5$ using different geometric configurations.}\label{fig:ex3_rec_1}
\end{figure}

In Figures \ref{fig:ex3_conf1_3d_65}, \ref{fig:ex3_conf2_3d_65}, \ref{fig:ex3_conf3_3d_65} and \ref{fig:ex3_conf4_3d_65}, we present the reconstructions at frequency $k=6.5$ for the configurations 1, 2, 3 and 4, respectively. In addition, in Figures \ref{fig:ex3_conf1_cs} 1, \ref{fig:ex3_conf2_cs} 2, \ref{fig:ex3_conf3_cs} 3 and \ref{fig:ex3_conf4_cs}, the cross section of the original obstacle and the reconstructions at $k=3.25$ and $k=6.5$ are shown for configurations 1, 2, 3 and 4, respectively.

The reconstruction obtained using configuration 1 is more accurate than the reconstructions obtained using the other configurations. This behavior was expected since the placement of the receptors in configuration 1 allows for obtaining information from a larger part of the obstacle. 

\subsection{Example 4: Using limited number of parameters} \label{ex:4}
This example is a continuation of Example 3. We use the same scattered data with 2\% noise that was generated for the configuration 1 in Example 3 to recover the obstacle in Figure \ref{fig:ex3_orig}. In Example 3, at each frequency, the number of modes used to approximate the boundary of the obstacle and the update obtained by the Gauss-Newton step is $N_p=\lfloor2k\rfloor$. Since the frequency varied from $k=0.5$ to $6.5$, the number of modes used varied from $1$ to $13$.

In Example 4, instead of letting the number of modes increase freely with the frequency, we set it to be $N_p=\min\left\{N_{max},\lfloor 2k \rfloor\right\}$, where we choose $N_{max}=8$, $10$ and $13$. All the other parameters for both the RLA and the Gauss-Newton method are the same as in Example 3.

\begin{figure}
\center
\begin{subfigure}{.45\textwidth}
\includegraphics[scale=0.2]{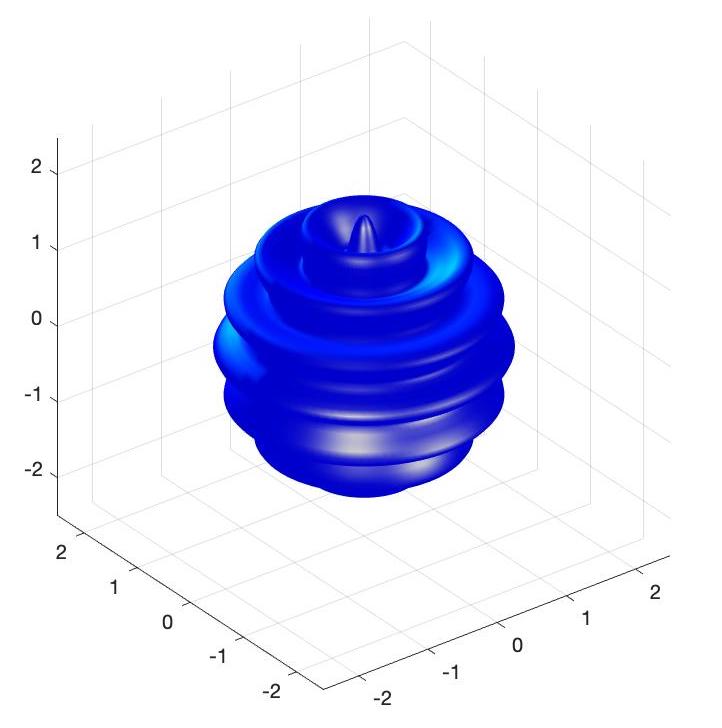}
\caption{8 modes}\label{fig:ex4_3d_h8}
\end{subfigure}
\begin{subfigure}{.45\textwidth}
\includegraphics[scale=0.2]{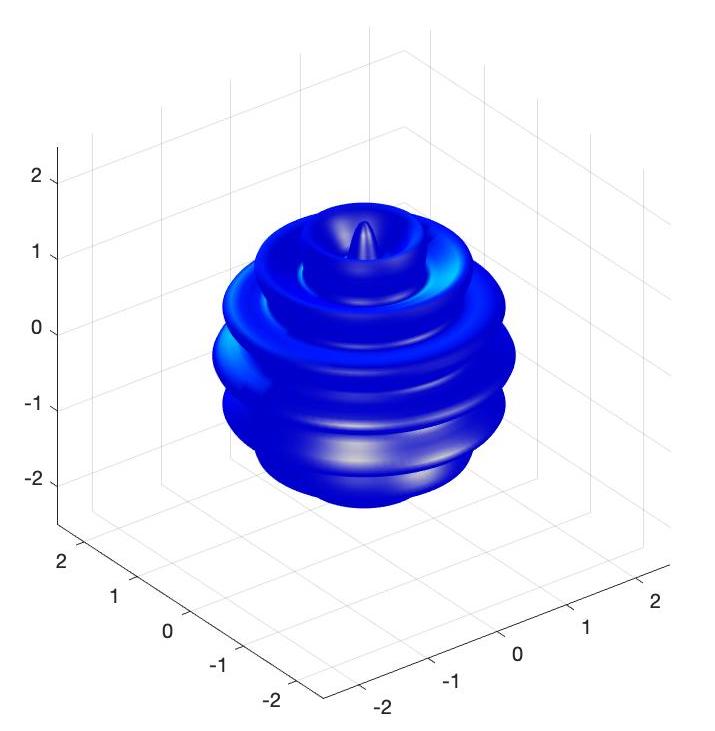}
\caption{10 modes}\label{fig:ex4_3d_h10}
\end{subfigure}

\begin{subfigure}{.45\textwidth}
\includegraphics[scale=0.2]{example_2a_iso.jpg}
\caption{13 modes}\label{fig:ex4_3d_h13}
\end{subfigure}
\begin{subfigure}{.45\textwidth}
\includegraphics[width=.8\textwidth]{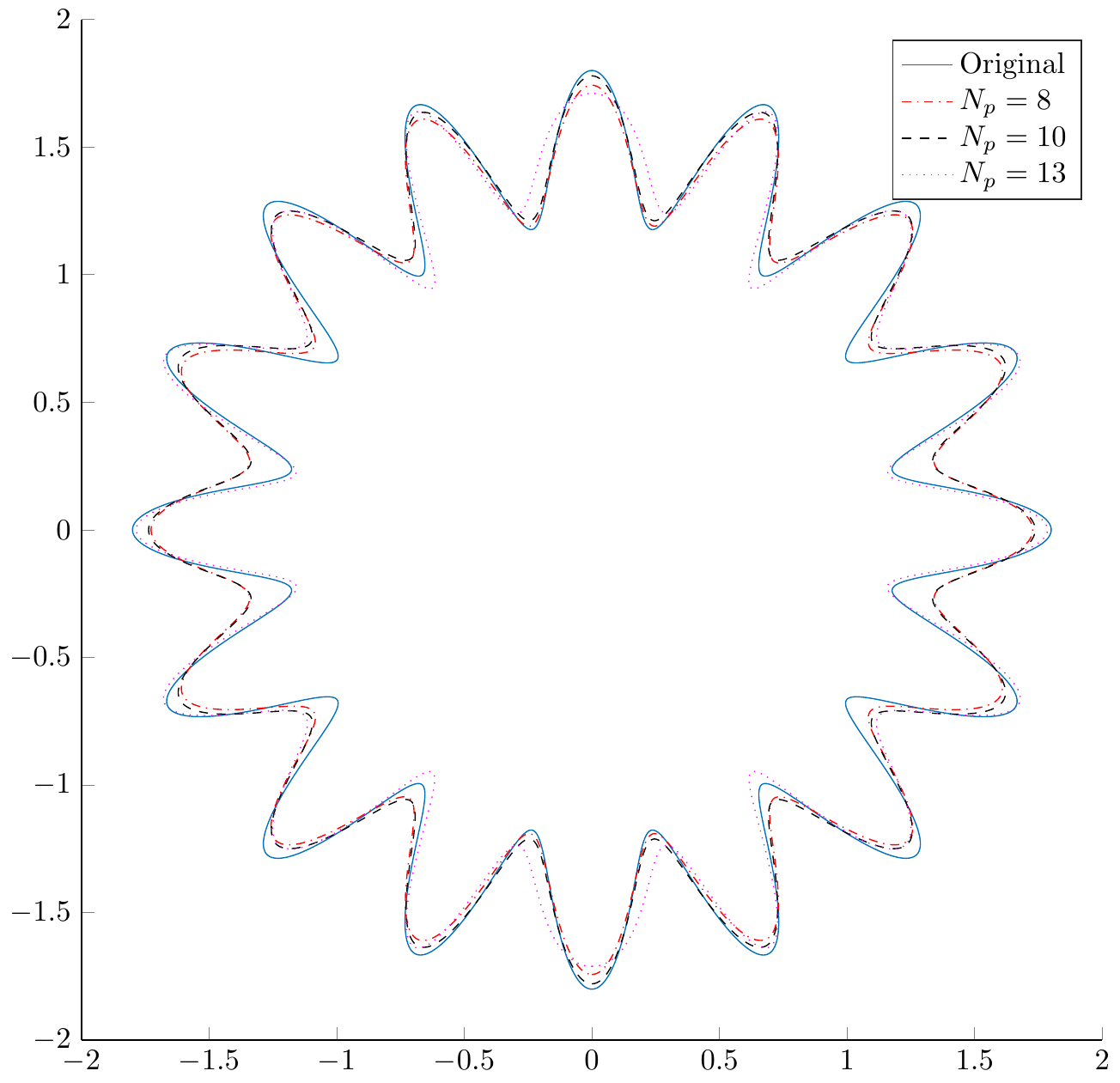}
\caption{Cross section}\label{fig:ex4_cs}
\end{subfigure}
\caption{{\bf (Example 4)} Reconstructions of the obstacle at wavenumber $k=6.5$ using the maximum number of modes equal to: (a) 8, (b) 10 and (c) 13. In (d), we present the cross section of the reconstructions for all modes and the original obstacle.}\label{fig:ex4_rec}
\end{figure}

In Figures \ref{fig:ex4_3d_h8}, \ref{fig:ex4_3d_h10} and \ref{fig:ex4_3d_h13}, we present the reconstructions at frequency $k=6.5$ using $8$, $10$ and $13$ modes, respectively. Figure \ref{fig:ex4_cs} has the cross section of the original obstacle and the reconstructions using $8$, $10$ and $13$ modes. 

As expected the reconstruction using 8 modes is more precise than the other reconstructions.  As we increase the number of modes used for the reconstruction, the results become increasingly worse due to the oscillations introduced by the higher order modes. An appropriate filter is required to damp the oscillation, which will be illustrated in the next example. 

\subsection{Example 5: Reconstruction of sharp features using multiple frequencies} \label{ex:5}
\begin{figure}
\center
\begin{subfigure}{.4\textwidth}
\includegraphics[scale=0.16]{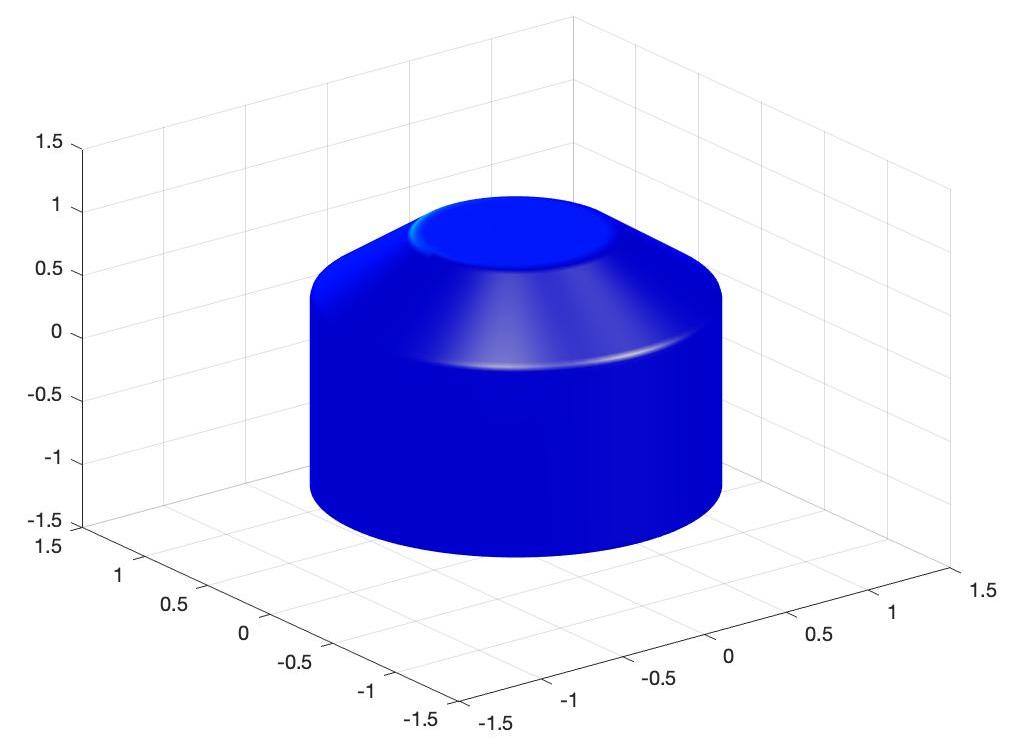}
\caption{Mine}\label{fig:ex5_orig}
\end{subfigure}
\begin{subfigure}{.4\textwidth}
\includegraphics[scale=0.09]{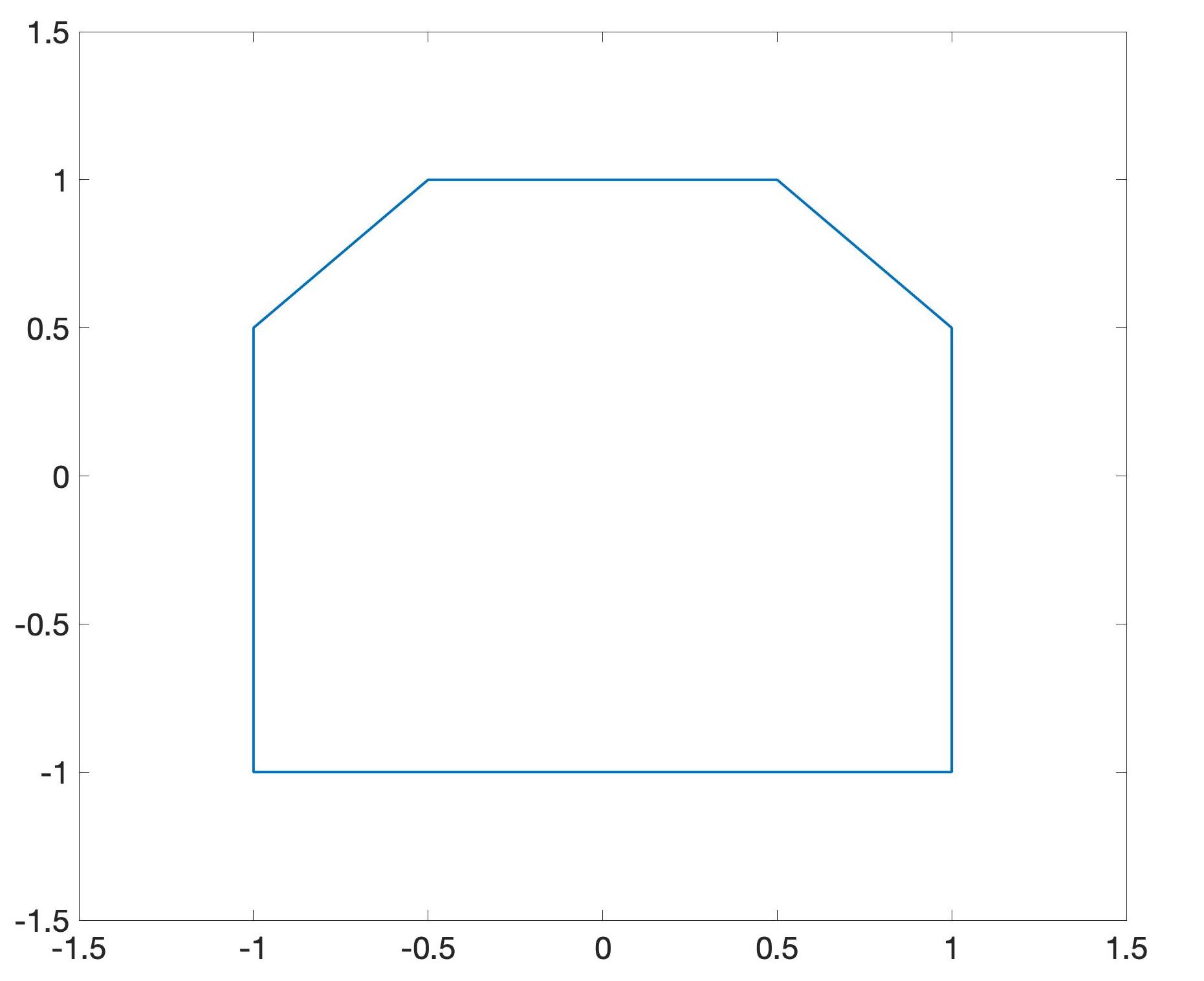}
\caption{Cross section of the mine}\label{fig:ex5_cs}
\end{subfigure}
\caption{{\bf (Example 5)} We present in (a) the original 3D obstacle and (b) a cross section of the obstacle.}\label{fig:ex5_general}
\end{figure}

In this example, we use multifrequency scattered data to reconstruct an obstacle with the shape of a land mine, see Figure \ref{fig:ex5_general}. One must use a vary large number of modes $N_p$ to recover the sharp corners of the obstacle using a trigonometric representation.
        
To generate the simulated scattered data ${\bf u}_{k_j,\db_{s}}^\emph{meas}$, we used incident plane waves given by $u_{k_j,\db_{s}}^\emph{inc}(\xb)$ $=e^{i k_j \db_{s} \cdot\xb}$, where $k_j=0.5+(j-1)0.25$, $j=1, \dots,78$, and
\begin{equation*}
\db_{s}=\left(\sin\left( m\pi/6 \right)\cos\left( 2n\pi/5 \right),\sin\left( m\pi/6 \right)\sin\left(  2n\pi/5 \right),\cos\left( m\pi/6 \right)\right),
\end{equation*}
with $s=(m-1)5+n$,  $m,n=1,\ldots,5$ and beyond that $\db_{26}=\left(0,0-1\right)$, and $\db_{27}=(0,0,1)$. The scattered field is measured at $N_r=900$ receptors located at the points 
\begin{equation*}
\xb_{lq}=10\left(\sin\left( l\pi/31\right)\cos\left( 2q\pi/30 \right),\sin\left( l\pi/31 \right)\sin\left( 2q\pi/30 \right),\cos( \pi/31)\right),
\end{equation*} 
with $l,q=1,\ldots,30$. As in Examples 3 and 4, to avoid inverse crimes, we add 2\% noise to the scattered data using formula \ref{eq:scat_noise}.

We apply the RLA with the damped Gauss-Newton method at each frequency. We used the same stopping criteria for the Gauss-Newton method as in Examples 3 and 4. We set the damping parameter $\alpha=0.1$ for the damped Gauss-Newton method at the initial frequency $k=0.5$, and $\alpha=0.1/(k\|h\|)$ for all other frequencies, where $\|h\|$ is the 2-norm of the vector of coefficients of the update. Regarding the number of modes in the polynomial representing the update $h$, we set it to be $N_p=\lfloor2 k\rfloor$. 

As we apply the RLA and increase the frequency, we note that oscillations are introduced in the reconstruction. These oscillations are an effect of the Gibbs phenomenon. They occur due to the limited number of modes used to recover the sharp edges of the obstacle. To mitigate the effect of the oscillations, we introduce an extra step in our reconstruction algorithm. After finding the domain update step $h$ for the damped Gauss-Newton method, we apply a low-pass Gaussian filter as follows:
\begin{equation*}
h(t)=h^{c}_0+\sum_{j=1}^{N_p}\ell_j \left(h^{c}_{j}\cos\left(2\pi j(t-0.5)\right)+h^{s}_{j}\sin\left(2\pi j(t-0.5)\right) \right) 
\end{equation*}
where the filter constants are $\ell_j=\exp\left(-(j/N_p)^2/\sigma^2\right)$, $j=1, \ldots, N_p$ and $\sigma^2$ is a constant to define the damping of the filter. The application of this low-pass Gaussian filter follows a similar logic as in \cite{Borges2015}, which uses a curve smoother developed in \cite{beylkin2014fitting}.

\begin{figure}[!htbp]
\center
\begin{subfigure}{.4\textwidth}
\includegraphics[scale=0.07]{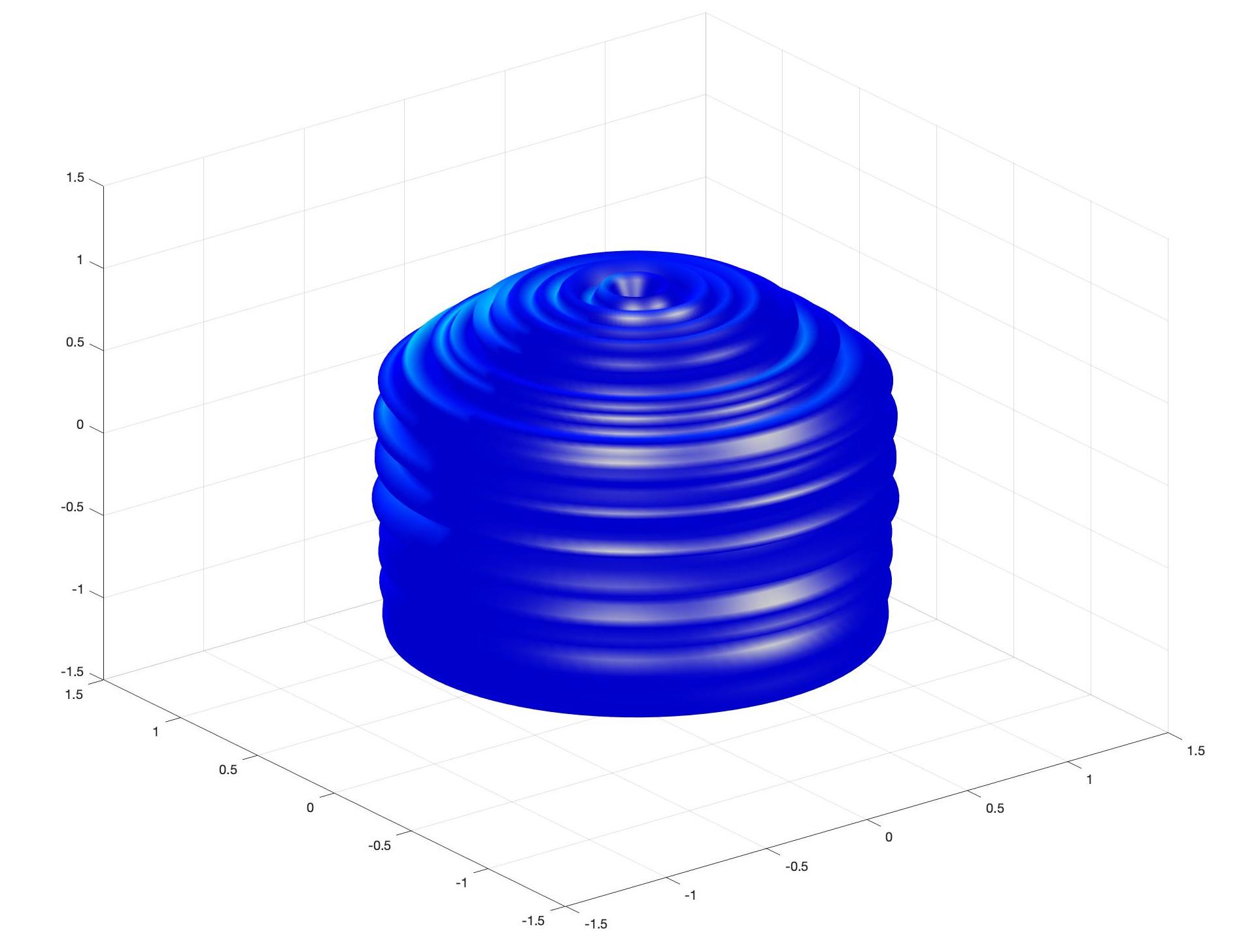}
\caption{No filter}\label{fig:ex5_3d_nofilter}
\end{subfigure}
\begin{subfigure}{.4\textwidth}
\includegraphics[scale=0.31]{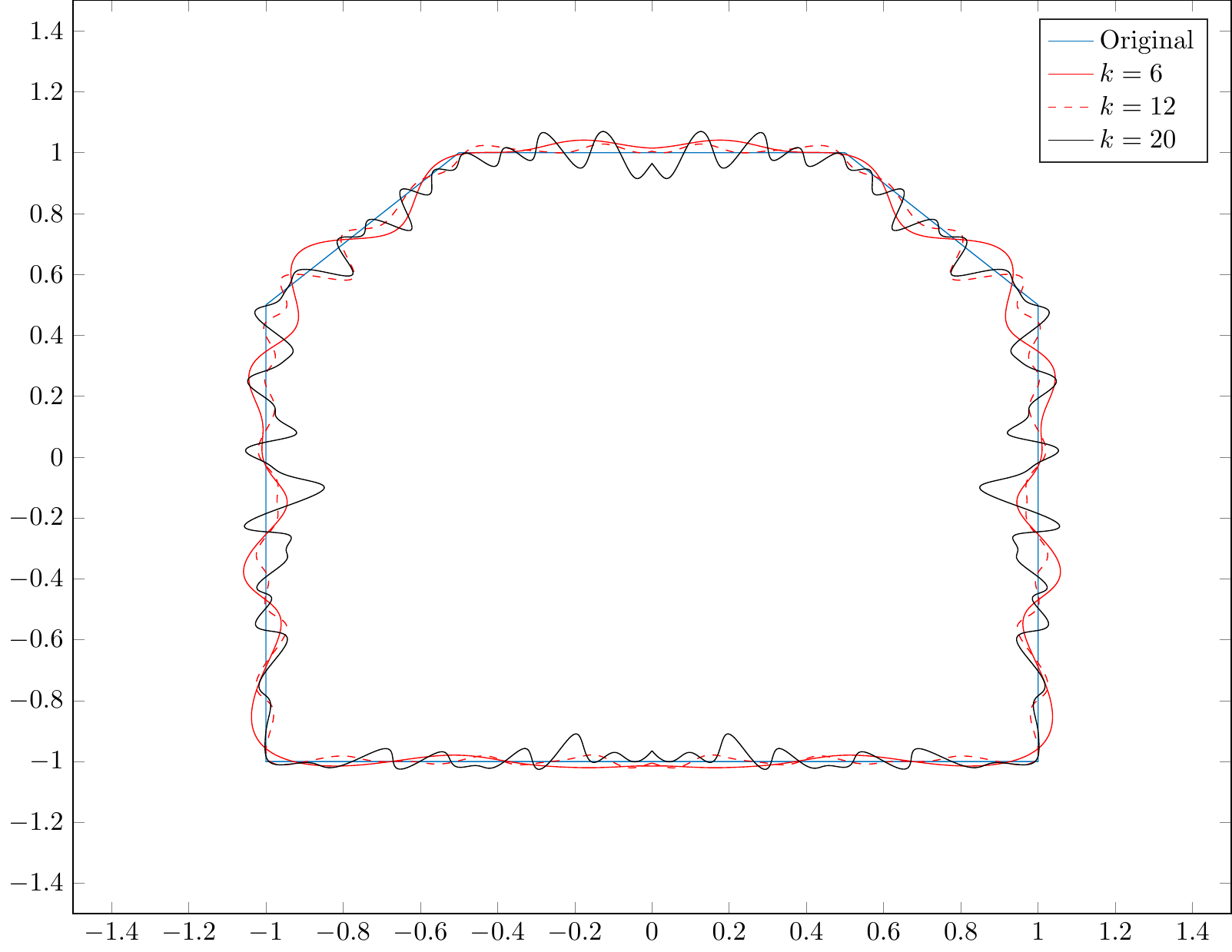}
\caption{No filter}\label{fig:ex5_cs_nofilter}
\end{subfigure}

\begin{subfigure}{.4\textwidth}
\includegraphics[scale=0.08]{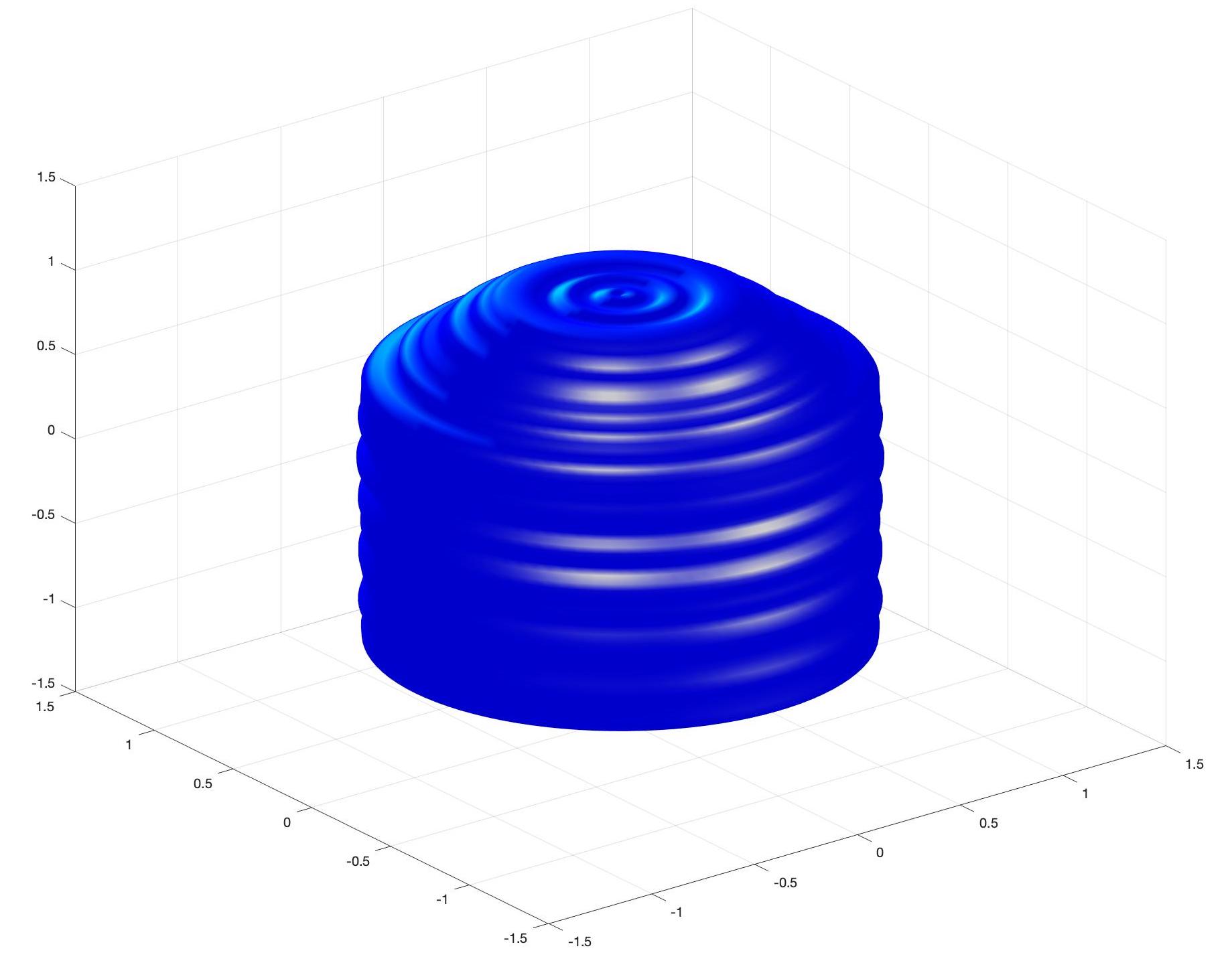}
\caption{Filter with $\sigma^2=0.5$}\label{fig:ex5_3d_sigma5}
\end{subfigure}
\begin{subfigure}{.4\textwidth}
\includegraphics[scale=0.33]{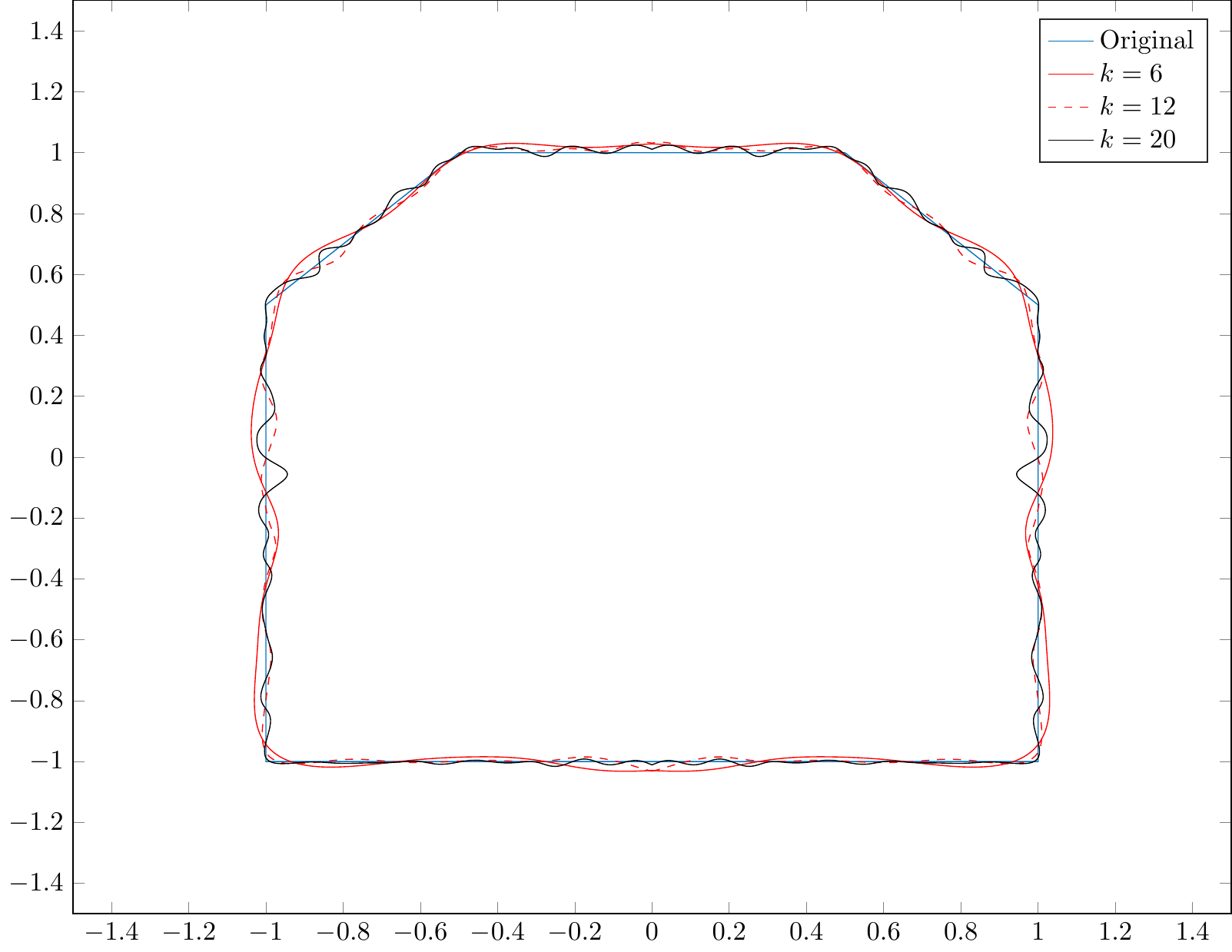}
\caption{Filter with $\sigma^2=0.5$}\label{fig:ex5_cs_sigma5}
\end{subfigure}

\begin{subfigure}{.4\textwidth}
\includegraphics[scale=0.08]{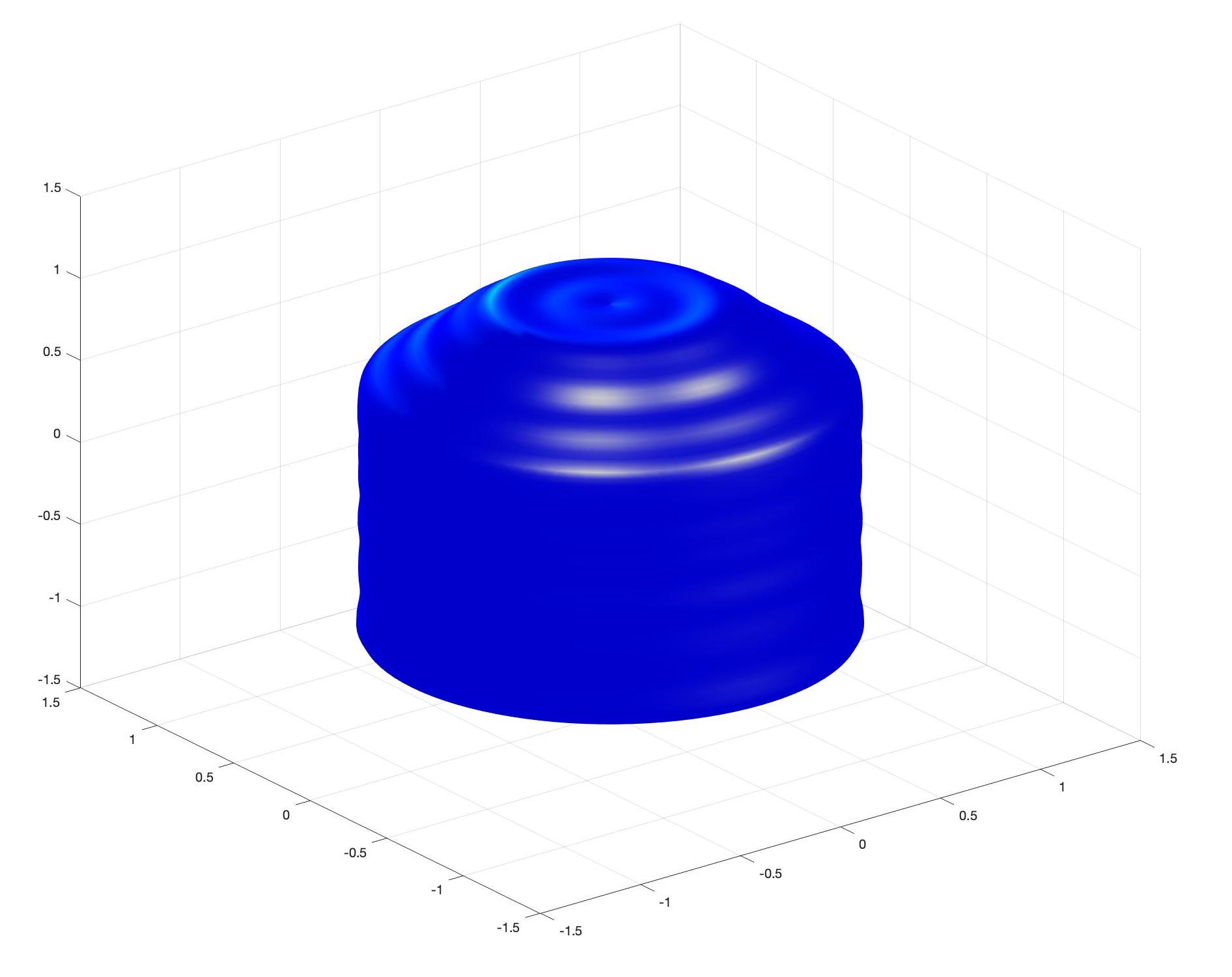}
\caption{Filter with $\sigma^2=0.1$}\label{fig:ex5_3d_sigma1}
\end{subfigure}
\begin{subfigure}{.4\textwidth}
\includegraphics[scale=0.33]{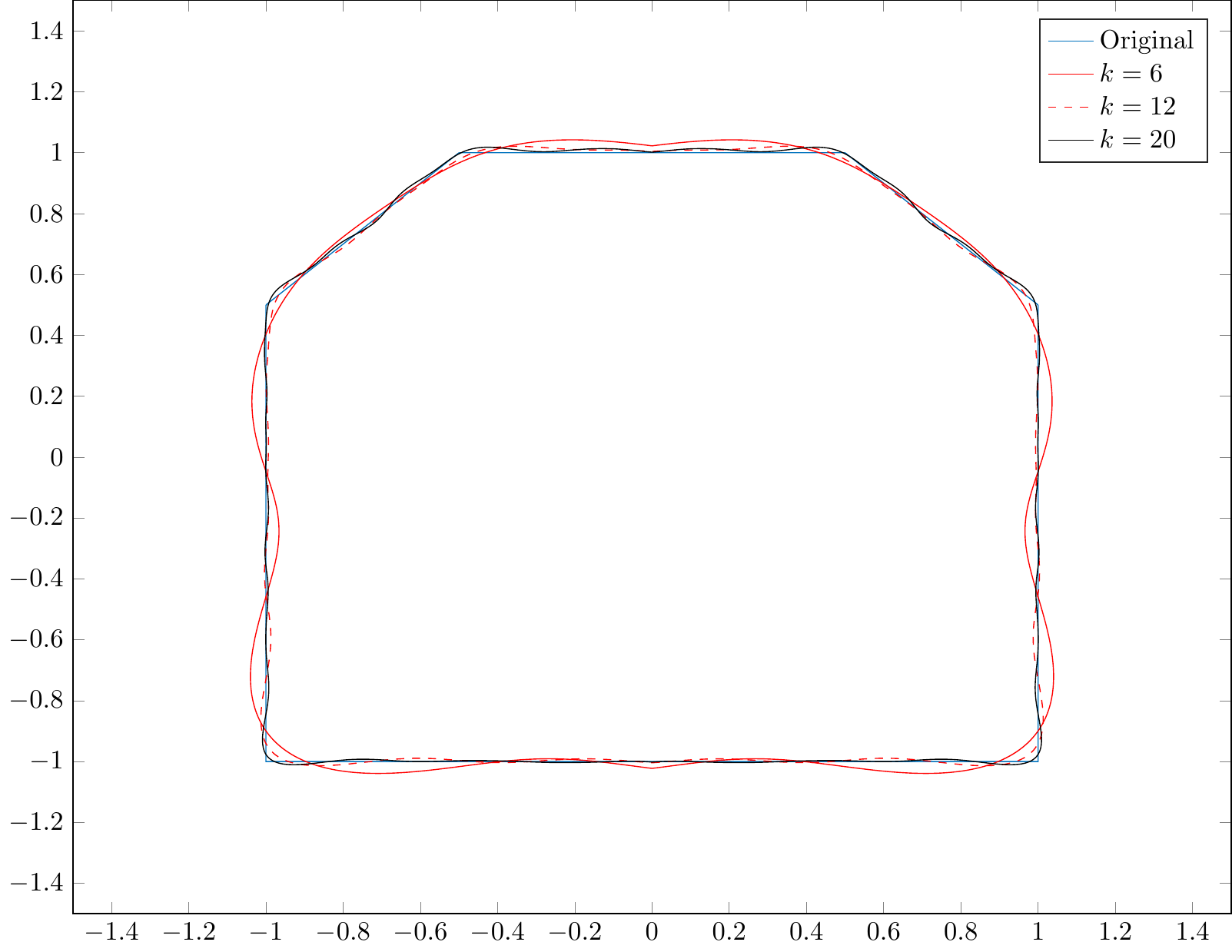}
\caption{Filter with $\sigma^2=0.1$}\label{fig:ex5_cs_sigma1}
\end{subfigure}
\caption{{\bf (Example 5)} Reconstructions of the obstacle at wavenumber $k=20$ are presented for the case when we use: (a) no filter, (c) filter with $\sigma^2=0.5$ and (e) filter with $\sigma^2=0.1$. The cross section of the original obstacle and the reconstructions at wavenumbers $k=6$, $12$ and $20$ for the cases when we use (b) no filter, (d) filter with $\sigma^2=0.5$ and (f) filter with $\sigma^2=0.1$.}\label{fig:ex5_rec}
\end{figure}

In Figures \ref{fig:ex5_3d_nofilter}, \ref{fig:ex5_3d_sigma5}, and \ref{fig:ex5_3d_sigma1},  we present the reconstruction obtained at $k=20$ using no filter, filter with $\sigma^2=0.5$, and filter with $\sigma^2=0.1$, respectively. In Figures \ref{fig:ex5_cs_nofilter}, \ref{fig:ex5_cs_sigma5}, and \ref{fig:ex5_cs_sigma1} we present the cross section of the original obstacle and of the reconstructions at $k=6$, $k=12$ and $k=20$ using no filter, filter with $\sigma^2=0.5$, and filter with $\sigma^2=0.1$, respectively.

The results clearly show that although the general shape of the obstacle can be recovered without filtering, it is hard to recover the sharp features of the obstacle with high resolution. On the other hand, using an appropriate filter produces a sharp reconstruction of the obstacle with very few oscillations.

\section{Conclusions}\label{s:conclusions}
In this paper, the forward and inverse scattering problems of recovering the shape of a three-dimensional impenetrable axis-symmetric sound-soft obstacle are studied. To solve the forward problem,  we make use of the symmetry of the obstacle by applying separation of variables in the azimuthal angle and Fourier decomposing the resulting problem. The original integral equation becomes a sequence of uncoupled line integral equations, where the new problems are both simpler and computationally cheaper to solve. For the inverse problem, we introduce a two-part framework for recovering the shape of the obstacle. In part 1, we find the axis of symmetry and center of the obstacle using the symmetry of the far field pattern. In part 2, we apply the RLA to obtain a reconstruction based on multifrequency data. 

We present five examples to examine the feasibility of the two-part framework. In Example 1, part 1 of the framework is tested successfully to obtain the axis of symmetry of an oblique ellipsoid. In Example 2, we show the interplay between the frequency and the local sets of convexity of the objective functional when the object is a sphere. In Example 3, we study different geometric configurations concerning the location of the receptors and the direction of the incident wave. In Example 4, we show that results are improved when the correct number of modes is used to represent the solution. Finally, in Example 5, we reconstruct an object with sharp edges using the RLA. A filter is used in the update of the domain to obtain an oscillation free high resolution reconstruction of the obstacle.

In the future, we intend to expand the inverse problem techniques to solve the multifrequency inverse scattering problem for three dimensional obstacles of arbitrary shape.


\FloatBarrier

\bibliography{./Biblio}

\end{document}